\documentclass{article}

\usepackage{arxiv}
\usepackage{amssymb}
\usepackage{amsmath}
\usepackage{empheq}
\usepackage{amscd}
\usepackage{graphicx}
\usepackage{graphics}
\usepackage[noadjust]{cite}
\usepackage{amsthm}
\usepackage{caption}
\usepackage{multirow}
\usepackage{array}
\usepackage{rotating}

\usepackage{indentfirst}
\usepackage{tikz}
\usepackage{calc}
\usepackage{epsfig}
\usepackage{natbib}
\usepackage[makeroom]{cancel}
\usepackage{multicol}
\usepackage{csquotes}
\usepackage{algorithm}
\usepackage{algorithmic}
\usepackage{enumitem}
\usepackage{hyperref}
\usepackage{url}
\hypersetup{colorlinks,linkcolor={blue},citecolor={blue},urlcolor={blue}}

\usepackage{timet}
\usepackage{graphicx}
\usepackage{multirow}
\usepackage{multicol}
\usepackage{lipsum}
\usepackage{timet}
\usepackage{epsfig}
\usepackage{amsmath}
\usepackage{amsfonts}
\usepackage{amssymb}
\usepackage{color}
\numberwithin{equation}{section}
\usepackage{caption}
\usepackage{float}
\usepackage{subcaption}
\usepackage{graphics}
\usepackage{xcolor,graphicx}
\usepackage{csquotes}
\usepackage{mathtools}
\usepackage{optidef}

\newcommand{\vertiii}[1]{{\left\vert\kern-0.25ex\left\vert\kern-0.25ex\left\vert #1 
    \right\vert\kern-0.25ex\right\vert\kern-0.25ex\right\vert}}

\usepackage{hyperref}       

\begin{document}
\title{Three-dimensional finite-difference \& finite-element frequency-domain wave simulation with multi-level optimized additive Schwarz domain-decomposition preconditioner: A tool for FWI of sparse node datasets}

%

\author{
  \href{https://orcid.org/0000-0002-6830-3767}{\includegraphics[scale=0.06]{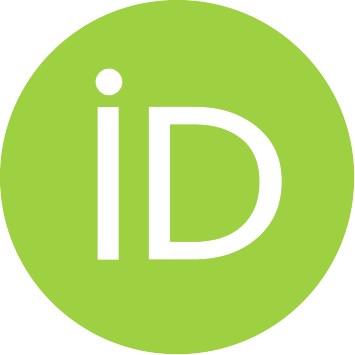}\hspace{1mm}P.-H. Tournier} \\
  Sorbonne University, CNRS, France.
  \texttt{tournier@ljll.math.upmc.fr} \\ 
  \And
  \href{https://orcid.org/}{\includegraphics[scale=0.06]{orcid.pdf}\hspace{1mm}P. Jolivet} \\
  IRIT-CNRS, University of Toulouse, France.
  \texttt{Pierre.Jolivet@enseeiht.fr} \\ 
  \And
  \href{https://orcid.org/}{\includegraphics[scale=0.06]{orcid.pdf}\hspace{1mm}V. Dolean} \\
  Univ. of Strathclyde, United Kingdom/UCA, CNRS, LJAD, France.
  \texttt{Victorita.Dolean@univ-cotedazur.fr} \\ 
  \And
\href{http://orcid.org/0000-0003-1805-1132}{\includegraphics[scale=0.06]{orcid.pdf}\hspace{1mm}H. S. Aghamiry} \\
  University Cote d'Azur - CNRS - IRD - OCA, Geoazur, Valbonne, France. 
  \texttt{aghamiry@geoazur.unice.fr}
\And
\href{http://orcid.org/0000-0002-4981-4967}{\includegraphics[scale=0.06]{orcid.pdf}\hspace{1mm}S. Operto} \\ 
  University Cote d'Azur - CNRS - IRD - OCA, Geoazur, Valbonne, France. 
  \texttt{operto@geoazur.unice.fr}
\And
\href{http://orcid.org/}{\includegraphics[scale=0.06]{orcid.pdf}\hspace{1mm}S. Riffo} \\ 
  University Cote d'Azur - CNRS - IRD - OCA, Geoazur, Valbonne, France. 
  \texttt{riffo@geoazur.unice.fr} 
  }

\renewcommand{\shorttitle}{Frequency-domain seismic wave solvers, Tournier et al.}

\maketitle

\begin{abstract}
Efficient frequency-domain Full Waveform Inversion (FWI) of long-offset node data can be designed with a few discrete frequencies hence allowing for compact volume of data to be managed. Moreover, attenuation effects can be straightforwardly implemented in the forward problem without computational overhead. However, 3D frequency-domain seismic modeling is challenging since it requires solving a large and sparse linear indefinite system per frequency with multiple right-hand sides. This linear system can be solved by direct or iterative methods. The former are very efficient to process multiple right-hand sides but may suffer from limited scalability for very large problems. Iterative methods equipped with a domain decomposition preconditioner provide a suitable alternative to process large computational domains for sparse node acquisition. 
The domain decomposition preconditioner relies on an optimized restricted additive Schwarz (ORAS) method, where a Robin or Perfectly-Matched Layer (PML) condition is implemented at the boundaries between the subdomains. The preconditioned system is solved by a Krylov subspace method, while a block low-rank Lower-Upper (LU) decomposition of the local matrices is performed at a preprocessing stage.  Multiple sources are processed in group with a pseudo-block method.
The accuracy, the computational cost and the scalability of the ORAS solver are assessed against several realistic benchmarks. For the considered benchmarks, a compact wavelength-adaptive 27-point finite-difference stencil on regular Cartesian grid shows better accuracy and computational efficiency than a P3 finite-element method on $h$-adaptive tetrahedral mesh, which remains however beneficial to comply with known boundaries such as bathymetry. The scalability of the method, the block processing of multiple RHS, the straightforward implementation of attenuation, which further improves the convergence of the iterative solver, make the method to be a versatile forward engine for large-scale 3D FWI applications from sparse node data sets.
\end{abstract}

\section{INTRODUCTION}
Frequency-domain full-waveform inversion (FWI) is suitable for long-offset stationary-recording acquisition, since subsurface models can be reconstructed with a few frequencies and attenuation is easily implemented without computational overhead \citep{Pratt_1996_TDV,Sirgue_2004_EWI,Gorszczyk_2017_TRW}.
Frequency-domain FWI is a data-fitting problem for subsurface parameter estimation constrained by a Helmholtz-type boundary-value equation.
This elliptic partial-differential equation is commonly discretized with finite differences or through more elaborate finite elements, which are crafted according to a specific problem \citep{Marfurt_1984_AFF}. Its discretization leads to a sparse linear system per frequency
whose unknown is the wavefield, the right-hand side (RHS) the seismic source and the matrix coefficients embed the subsurface properties. To solve such a system, we can use either a sparse direct solver \citep{Duff:2017:DMS} or an iterative one \citep{Saad:2003:IMS}.
While a direct solver is efficient when processing multiple RHSs for problems involving less than $\sim 100$ million unknowns \citep{Amestoy_2016_FFF,Kostin_2019_DFA,Amestoy_2021_EDA}, the memory overhead generated by the storage of the LU factors and the limited scalability of the LU decomposition make application on large scale problems difficult. The second approach relies on iterative solvers, whose natural scalability and moderate memory demand make them suitable for large-scale problems (see \citet{Plessix_2017_CAT} for a comprehensive discussion). 

However, due to the indefiniteness of the operator, fast convergence of iterative solvers for the discretized Helmholtz equation remains challenging as shown in \citet{Ernst:12:NAM}, hence requiring efficient preconditioners. 
To tackle this issue, many techniques have been developed, which can be divided into three broad categories: incomplete factorizations, e.g., analytic incomplete LU \citep{Gander:2000:AILU}, multigrid with shifted Laplacian \citep{Riyanti_2007_PMP,Erlangga:2008:AIM}, and domain decomposition preconditioners including overlapping Schwarz type methods among them \citep{Dolean:2015:IDD}. In turn, domain decomposition methods can be of different types. We refer the reader to the paper of \citet{Gander:2018:SIREV} for an overview of these solvers and the link between them. 

Domain decomposition offers a very suitable parallel framework to tackle large scale problems repeatedly, as in the case of inverse problems with multiple RHSs. In this context, an example of the use of a domain decomposition method as a solver for a medical imaging problem based on time-harmonic Maxwell's equations can be found in \citet{Tournier:2019:MTI} and due to the similarities between Maxwell's and acoustic wave equation \citep{Carcione_2002_AEA}, a similar type of approach can be applied to Helmholtz equations \citep{Ben-Hadj-Ali_2011_GEO}. 

In this paper, we present a new solver for 3D frequency-domain wave simulation as a forward engine for FWI, which consists in combining the Krylov subspace generalized minimal residual (GMRES) solver \citep{Saad:2003:IMS} with an optimized restricted additive Schwarz (ORAS) domain-decomposition preconditioner \citep{Graham:2017:RRD,Bonazzoli:2019:ADD}. We will present both one-level and two-level methods and discuss their domain of application.
Two-level methods are obtained from one-level methods by adding a coarse space correction in order to achieve scalability with respect to the number of subdomains. However, a specificity of wave propagation problems, in particular those arising in exploration geophysics, is the sensitivity of the method to the physical parameters such as the range of wave numbers and the level of heterogeneity of the medium \citep{Plessix_2017_CAT}. To manage this issue, several coarse space correction strategies have been investigated for Helmholtz problems including those of interest here in \citet{Bootland:2021:ACS,Bootland:2021:SWA} and the coarse space based on a geometrical grid (aka the grid coarse space) was considered to be the best strategy for solving the applications of interest especially in the case of parsimonious discretization with the lowest number of points per wavelength possible. This choice was pushed further and assessed in more detail on a few test cases in \citet{Dolean:2020:IFD,Dolean:2020:LSF}. On the other side, spectral coarse spaces like those based upon Dirichlet-to-Neumann (DtN) operator \citep{Conen:2014:ACS,Bootland:2019:ODN} remain potentially promising alternatives for inverse problems involving multiple RHSs where the parallel precomputation of the coarse space components could bring a clear advantage and improve considerably the overall computation time. 
Compared to the well-known preconditioner based upon shifted Laplacian \citep{Erlangga:2008:AIM}, let's add that the optimized Schwarz preconditioners we introduce here are less sensitive to the shift (added attenuation) and can be used without it. 

In this study, we assess two discretization schemes of the Helmholtz equation. The first relies on continuous Lagrange finite elements of order 3 (P3) on unstructured tetrahedral meshes to comply with complex geometries and adapt the size of the elements to the local wavelength ($h$-adaptivity). Our rationale behind the choice of the polynomial order is to choose the coarsest discretization providing the desired accuracy while being able to capture heterogeneities during FWI whose size is of the order of half the wavelength, namely the theoretical resolution of FWI \citep{Virieux_2009_OFW}.
The second discretization scheme relies on the recently proposed compact wavelength-adaptive 27-point finite-difference method on uniform Cartesian grid \citep{Aghamiry_2021_AFD}. The governing idea of this method compared to the classical 27-point stencil \citep{Operto_2007_FDFD} is to make its accuracy uniform with respect to the local wavelength through a local adaptation of the stencil coefficients. Finite difference methods are the most popular ones in exploration geophysics due to their simplicity and their ability to perform well on structured meshes with the lowest number of degrees of freedom for a prescribed accuracy. On the other side, finite elements may be favored when their versatility for dealing with unstructured geometries and non-trivial topographies is needed.  The different paradigms governing these two methods prompt us to provide a re-assessment of their pros and cons in the framework of Helmholtz problems.

This paper is organized as follows. In the first section, we review the ORAS preconditioner and its two-level variant with the grid coarse space.  We present their formulation and review several key algorithmic aspects. Then, we briefly review the two discretization schemes and their dispersion properties. The third section presents numerical results. We first define the optimal tuning of the ORAS solver that achieves the best compromise between accuracy and computational efficiency. Then, we assess the two discretization schemes against four large-scale benchmarks. The solutions are validated against analytical solutions for the first two benchmarks and the convergent Born series (CBS) method \citep{Osnabrugge_2016_CBS} for the last two benchmarks, which involve complex media. We conclude the numerical section with a weak and strong scalability analysis. In the final discussion section, we discuss the pros and cons of the proposed method against alternative approaches based on sparse direct solver and explicit time marching methods with a numerical example.


\section{Method}
We consider the Helmholtz equation, which is the simplest mathematical model of frequency-domain acoustic wave propagation:
\begin{equation}
\label{eq1}
\left(\Delta + k^2(\bold{x}) \right) u(\bold{x},\omega) = f(\bold{x},\omega), ~ \text{in a subsurface domain } \Omega,
\end{equation}
where $u$ is the monochromatic pressure wavefield, $f$ the monochromatic source, $k(\bold{x},\omega)=\omega/c(\bold{x})$ is the wavenumber with $\omega$ denoting the angular frequency, $c(\bold{x})$ the wavespeed (which is complex valued in viscous media) and $\bold{x}=(x,y,z) \in \Omega$. 
Let us denote its boundary by $\partial\Omega = \Gamma_{D} \cup \Gamma_{R}$ where Dirichlet (free-surface) conditions are imposed on $\Gamma_{D}$ ($u = 0$) and a Robin condition on $\Gamma_{R}$ ($\frac{\partial u}{\partial \boldsymbol{n}} + i k u  = 0$). The Robin condition is a standard first order approximation to the far field Sommerfeld radiation condition enabling the description of the wave behavior in a bounded domain, while the physical domain is not bounded. Alternatively, perfectly-matched layer boundary conditions can be applied on $\Gamma_{R}$ \citep{Berenger_1994_PML,Chew_1994_PMM}. The resulting linear system, while being complex symmetric, is not Hermitian and this will condition the choice of iterative method to solve it.  
The method that is reviewed below can be extended to more general form of the wave equation including heterogeneous density \citep[e.g.,][]{Operto_2007_FDFD}, anisotropic effects \citep[e.g.,][]{Operto_2014_FAT}, and elastic effects \citep[][]{Gosselin_2014_FDF,Li_2015_FDE,Li_2020_FEW}, with however a potential impact of the physics on the convergence behavior of the method \citep{Plessix_2017_CAT}.
In the sequel of this study, the discretized form of equation~\ref{eq1} is written in matrix form as
\begin{equation}
\label{eq1.1}
A\mathbf{u}=\mathbf{f}.
\end{equation}
where $\mathbf{u}$ is the vector of degrees of freedom.
\subsection{Optimized Schwarz preconditioners}\label{sec:DDM}

We review the optimized Schwarz preconditioners that we use to solve efficiently the linear system \eqref{eq1.1}. A well-known iterative solver for this type of indefinite linear systems is the Krylov subspace Generalized Minimal RESidual Method (GMRES) \citep{Saad:2003:IMS}. However, the Helmholtz operator requires efficient preconditioning, which can be done by domain decomposition \citep{Dolean:2015:IDD}. Note that we will employ right preconditioning to solve~\eqref{eq1.1}: we will solve 
\begin{equation}
A M^{-1} \mathbf{y}=\mathbf{f} \quad \text{where} \; \mathbf{u} = M^{-1} \mathbf{y}
\end{equation}
with preconditioner $M^{-1}$. Compared to the left-preconditioned variant, right-preconditioned GMRES is usually considered more stable, and the residual of the preconditioned system is independent of the preconditioner.

Our approach is based on the optimized restricted additive Schwarz (ORAS) method. Let us first introduce the one-level variant. In order to use the domain decomposition method, we first partition $\Omega$ into overlapping subdomains $\left\lbrace\Omega_{j}\right\rbrace_{j=1}^{N}$. Note that graph partitioners such as Metis \citep{Karypis_2013_MSP} and Scotch \citep{Francois_2018_SCOTCH} provide a non-overlapping partition from which the overlapping counterpart can be inferred by adding more layers in a recursive manner when subdomains with larger overlap are required. Now that we have a domain decomposition, we can define the restriction to a given subdomain $\Omega_{j}$ via a Boolean matrix $R_{j} \in \mathbb{R}^{n_{j} \times n}$ where $n_{j}$ is the number of degrees of freedom in $\Omega_{j}$. Since our subdomains overlap, we also make use of a partition of unity having matrix form $D_{j} \in \mathbb{R}^{n_{j} \times n_{j}}$, which is diagonal and satisfies $\sum_{j=1}^{N} R_{j}^{T} D_{j} R_{j} = I$. Note that $R_{j}^{T}$ acts as an extension by zero outside of $\Omega_{j}$. We can now define the one-level optimized restricted additive Schwarz (ORAS) preconditioner as 
$$
M_{\text{ORAS}}^{-1} = \sum_{j=1}^{N} R_{j}^{T} D_{j} B_{j}^{-1} R_{j},
$$
where $B_{j}$ is the local matrix on $\Omega_{j}$ and Robin or impedance boundary conditions are assumed on $\partial\Omega_{j} \setminus \partial\Omega$ (as mentioned previously, perfectly-matched layer boundary conditions can be used instead in the case of regular partitioning). Each contribution in the sum of the above equation can be computed locally in parallel. The ORAS method is a one-level method as it relies on local subdomain solutions and local transmission of data. This is reflected in the scalability of such a method, as the number of iterations usually grows as we increase the number of subdomains $N$. Therefore, in order to achieve robustness with respect to $N$, a coarse space needs to be included, leading to a two-level method. 

In a two-level method, the typical characteristics of the method (like the diameter $H$ of the elements in the coarse grid or the number of eigenvectors in the coarse space) can be tuned in order to achieve the desired robustness, i.e., the weak dependence on the wave number $k$ or the heterogeneities at the same time as the scalability. The typical two-level preconditioner is mainly characterized by the thin and long matrix $Z$ and can be written in a generic algebraic way as follows
\begin{equation}
\label{eq:BNN}
M^{-1}_{2,ORAS} =  M^{-1}_{ORAS} P + \underbrace{Z E^{-1} Z^{\dagger}}_Q.
\end{equation}
Here $M_{\text{ORAS}}^{-1}$ is the one-level preconditioner, $Z$ is a rectangular matrix with full column rank, $E = Z^{\dagger}A Z$ is the so-called coarse grid matrix, $Q = Z E^{-1} Z^{\dagger}$ is the so-called coarse grid correction matrix, and $P = I - A Q$ is a projection matrix. The most natural coarse space is the one based on a coarser mesh or grid, the so-called grid coarse space. In the case of finite elements for example, let us consider $\mathcal{T}^H$, a coarse simplicial mesh of $\Omega$ with mesh diameter $H$, and $V^{H} \subset V$, the corresponding finite element space. Let $Z$ the interpolation matrix between the coarse and the fine finite element spaces. Then, in this case, $E = Z^{\dagger}A Z$ is really the stiffness matrix of the problem discretized on the coarse mesh.
Note that the inclusion of the coarse space can be done in several ways, the one indicated in \eqref{eq:BNN} is based on an adapted deflation technique by making use of the projector $P$ outside the coarse space for an increased efficiency.

As mentioned previously, the two discretization schemes considered in this paper are wavelength-adaptive 27-point finite differences and P3 finite elements. In the case of finite differences, the solution is defined by its values on the nodes of the Cartesian grid with a resolution of four points per wavelength ; this unfortunately means that a coarser grid discretization would lead to a poor representation of the solution, invalidating the benefits of a second level of preconditioning based on a coarse grid. Thus, the one-level ORAS preconditioner will be used for the solution of linear systems stemming from finite-difference discretization.

However, in the case of P3 finite elements, the much larger set of degrees of freedom allows for a relevant physical representation of the solution even on a coarser mesh with a resolution of two points per wavelength. Thus, contrary to the finite difference case, the second-level based on a coarse mesh discretization for finite elements still maintains a pertinent physical meaning and has a significant positive impact on convergence. Hence, the two-level variant will be used for finite element discretizations.

In the next section, we give more details on the practical aspects of the parallel implementation of the one- and two-level ORAS preconditioners.

\subsection{Algorithmic aspects and parallel implementation}\label{sec:Algo} 

\paragraph{One-level method} The first step is to decompose the computational domain into overlapping subdomains. For finite-difference discretizations, the regular grid is partitioned into overlapping cuboids. For finite-element discretizations on unstructured meshes, we rely on automatic graph partitioners such as Metis and Scotch, which produce a non-overlapping decomposition of the set of mesh elements while minimizing interfaces between subdomains and conserving good load-balancing. The overlapping decomposition is then obtained by simply adding successive layers of elements to reach the desired size of overlap. The choice of the width of the overlap region is discussed in more details below.

The implementation of the domain decomposition method used here relies on Message Passing Interface (MPI) parallelism. Each subdomain of the decomposition is assigned to one MPI process (the number of subdomains $N$ is then equal to the number of MPI processes).

The parallel implementation operates only on distributed quantities; for instance, a global vector $\mathbf{v}$ is represented as a distributed vector, i.e., a collection of local vectors $(\mathbf{v}_j)_{1 \leq j \leq N}$ where the values on the unknowns shared between multiple subdomains are the same:
$$\mathbf{v}_j = R_j \mathbf{v}, \quad j = 1,\dots,N.$$

The MPI process in charge of subdomain $j$ then holds the restriction $R_j \mathbf{v}$ to subdomain $j$. Distributed operators can be defined in a similar manner; for example, the global operator $A$ can be written as
\begin{align}
    A =& \left( \sum_{i=1}^N R_i^T D_i R_i \right) A \\
        =& \sum_{i=1}^N R_i^T D_i A_i R_i,
\end{align}
where $A_i = R_i A R_i^T$ is the local restriction of $A$ to subdomain $i$ and the last equality holds under the hypothesis that $D_i$ vanishes at the subdomain interface ; we then have $D_i R_i A = D_i A_i R_i$.\\

Thus, the distributed result $(\mathbf{u}_j)_{1 \leq j \leq N}$ of the parallel matrix-vector product with $A$ applied to vector $\mathbf{v}$ distributed as $(\mathbf{v}_i)_{1 \leq i \leq N}$ is
\begin{align}
    \mathbf{u}_j =& R_j A \mathbf{v}\\
        =& R_j \sum_{i=1}^N R_i^T D_i A_i R_i \mathbf{v}\\
        =& R_j \sum_{i=1}^N R_i^T D_i A_i \mathbf{v}_i \\
        =& D_j A_j \mathbf{v}_j + \sum_{i \in \mathcal{O}(j)} R_j R_i^T D_i A_i \mathbf{v}_i,
\end{align}
where $\mathcal{O}(j)$ is the set of neighbors of subdomain $j$. Exchange operators $R_j R_i^T$ correspond to neighbor-to-neighbor MPI communications, transferring the values of the local input vector on the shared overlap region of subdomain $i$ to corresponding neighboring subdomains $j$. Thus, applying operator $A$ in parallel amounts to performing local weighted matrix-vector products $D_i A_i \mathbf{v}_i$ on each MPI process, exchanging the result in the overlapping region with neighboring subdomains, and summing contributions.\\

In the same manner, the distributed result $(\mathbf{u}_j)_{1 \leq j \leq N}$ of the application of the one-level ORAS preconditioner $M_{\text{ORAS}}^{-1}$ to a distributed vector $(\mathbf{v}_i)_{1 \leq i \leq N}$ is
\begin{align}
    \mathbf{u}_j =& R_j M_{\text{ORAS}}^{-1} \mathbf{v}\\
    =& R_j \sum_{i=1}^N R_i^T D_i B_i^{-1} R_i \mathbf{v}\\
    =& R_j \sum_{i=1}^N R_i^T D_i B_i^{-1} \mathbf{v}_i\\
    =& D_j B_j^{-1} \mathbf{v}_j + \sum_{i \in \mathcal{O}(j)} R_j R_i^T D_i B_i^{-1} \mathbf{v}_i,
\end{align}
which amounts to performing local weighted solves $D_i B_i^{-1} \mathbf{v}_i$ on each MPI process, exchanging the result in the overlapping region with neighboring subdomains, and summing contributions. The idea is then to precompute local factorizations of matrices $B_i$ in each subdomain using a direct solver such as the MUltifontal Massively Parallel Solver (MUMPS) \citep{MUMPS_2021_MMP}. This can be viewed as the \textit{setup} cost of the preconditioner. Then, during the iterative solution of the linear system, local solves $B_i^{-1} \mathbf{v}_i$ occurring in the application of the preconditioner can be performed by forward-backward substitution using the local factorization of $B_i$ in each subdomain.\\

All operations involved in the iterative solution of the linear system can be distributed the same way, and we have all the ingredients to solve the problem with parallel distributed preconditioned GMRES.

\paragraph{Note on the size of overlap} For wave propagation problems, the idea of the Optimized Schwarz method is to construct $B_i$ as a restriction of the global operator $A$ on subdomain $i$ with an approximation of transparent boundary conditions at the subdomain interface. In the case of finite-differences discretizations, a good approximation of transparent boundary conditions can be implemented using PMLs on the interface boundaries of the cuboid-shaped subdomains. In this case, the idea is to add a PML in each direction in the overlap region, with the width of the overlap sufficiently large for the PML to be a good enough approximation of transparent boundary conditions. In practice, we add three grid points in each direction (with a corresponding overlap width of six grid elements) for most of the experiments in this paper as it is found to be the best compromise (note that increasing the overlap width increases the computational cost), although we observed that the optimal amount of overlap can increase with frequency.

\paragraph{Note on processing multiple right-hand sides} Frequency-domain FWI requires forward simulations with multiple sources, which means solving linear systems with multiple right-hand sides (RHSs). Instead of opting for the naive approach which consists in solving the linear systems successively for each RHS, we can use various ways to accelerate computations when processing multiple RHSs, from pseudo-block Krylov methods to block Krylov method with recycling such as the Block Generalized Conjugate Residual Method with Inner Orthogonalization with inexact breakdown and deflated restarting (BGCRO-DR) \citep{Giraud_2021_BGC}. The simplest ones are pseudo-block methods, where operations for each RHS are simply fused together in order to achieve higher arithmetic intensity and decrease the number of global synchronizations. In particular, for our one-level ORAS method, this means that we can also leverage the multi-RHS processing capabilities of the direct solver for local problems when computing the application of the ORAS preconditioner $M_{\text{ORAS}}^{-1}$ to a block of input vectors, using the factorization of $B_i$ to perform forward elimination and backward substitution simultaneously on each input vector in order to compute the action of $B_i^{-1}$ on the whole block.

\paragraph{Two-level method} As mentioned previously, the coarse space we consider here is based on a coarse grid discretization of the PDE. Since a finite-difference grid with less than four points per wavelength cannot carry a meaningful representation of the solution, we only use the two-level method in the finite element case, where the larger set of degrees of freedom of the P3 discretization allows for a relevant representation of the solution even on a coarser mesh resolution of two points per wavelength.

For the two-level method, two meshes of the same domain must be considered simultaneously. Here we will consider nested meshes, so that the fine mesh is generated as a uniform refinement in which edges of all tetrahedra are divided uniformly in $s$. In this work we choose $s = 2$, which means that the coarse mesh is constructed with a resolution level of two points per wavelength. The required resolution of four points per wavelength will then be achieved on the fine mesh obtained by splitting.

The number of unknowns in the coarse problem is then $s^d = 8$ times smaller than the fine problem. However, for large frequencies the coarse problem still reaches tens of millions of unknowns. Thus, instead of relying on an exact factorization of the coarse matrix $E$, we solve it iteratively using a preconditioned GMRES method. The preconditioner for the coarse problem is again a one-level ORAS preconditioner ; this can be seen as a two-level recursive DD method. This outer-inner strategy requires the use of flexible methods such as FGMRES~\citep{saad1993flexible}. Moreover, we observe that in practice an accurate solution of the coarse problem is not required: the best compromise is obtained by choosing a loose relative GMRES tolerance of $10^{-1}$ for the coarse problem.

Finally, in order to minimize communications between levels, a natural choice is to use the same subdomain partitioning for the coarse and fine levels. This means that each MPI process will be in charge of two (the fine and the coarse) discretizations of the same spatial subdomain, with fine subdomains being merely uniformly refined subdomains of the coarse grid. Thus, by setting the width of the overlap to 1 for the coarse mesh partitioning, the overlap at the fine level will then have width $s = 2$.


%
%
%

\subsection{Discretization of the Helmholtz problem}
\label{sec:discret}
We now review the two discretization methods that we consider in this study. The first relies on finite elements with Lagrange polynomials while the second one relies on finite differences. 

\subsubsection{Finite-element frequency-domain (FEFD) discretization}
The first discretization scheme relies on the finite element method \citep{Zienkiewicz_2005_FEM}, with piecewise Lagrange polynomials on a simplicial mesh $\mathcal{T}^{h}$ of $\Omega$, which has a characteristic element diameter $h$.  
The weak form of the problem is to find $u \in V$ such that $a(u,v)  = F(v),\, \forall \ v \in V_{0}$
where
\begin{align}
\label{WeakFormBilinearPart}
a(u,v) & = \int_{\Omega} \left( \nabla u \cdot \nabla \bar{v} - k^2 u \bar{v}\right) \, \mathrm{d}x + \int_{\Gamma_{R}} i k u \bar{v} \, \mathrm{d}s,
& F(v) & = \int_{\Omega} f \bar{v} \, \mathrm{d}x,
\end{align}
and $V_{0} = \left\lbrace u \in H^{1}(\Omega) \colon u = 0\text{ on } \Gamma_{D}\right\rbrace$.
Denoting the associated trial space $V^{h} \subset H^{1}(\Omega)$ and test space $V_{0}^{h} \subset V_{0}$, the discrete problem is to find $u_{h} \in V^{h}$ such that $a(u_{h},v_{h}) = F(v_{h}),\, \forall \ v_{h} \in V_{0}^{h}.$ Let $\left\lbrace\phi_{j}\right\rbrace_{j=1}^{n}$ be the nodal basis for $V_{0}$. Then we can rewrite the discrete weak formulation as a (complex) linear system
\begin{align}
\label{LinearSystem}
A\mathbf{u}= \mathbf{f},
\end{align}
where the coefficient matrix $A \in \mathbb{C}^{n \times n}$ and right-hand side vector $\mathbf{f} \in \mathbb{C}^{n}$ are given by $A_{i,j} = a(\phi_{j},\phi_{i})$ and $\mathrm{f}_{i} = F(\phi_{i})$. We then seek the solution $\mathbf{u} \in \mathbb{C}^{n}$ of the (in general) complex symmetric indefinite system \eqref{LinearSystem} to give $u_{h}(x) = \sum_{j=1}^{n} \mathrm{u}_{j} \phi_{j}(x)$.

The oscillatory nature of solutions to the Helmholtz equation requires a sufficiently fine mesh in order to get a good approximation to the problem. Moreover, when the wave number $k$ is increasing, for the same level of accuracy the number of grid points must increase faster than $k$ increases, in order to avoid pollution effect \citep{Babuska:1997:IPE}. This growth depends on the discretization and becomes less important when the precision of the method increases in the case of finite element methods. For example, for piecewise linear (P1) finite element approximation on simplicial elements of diameter $h$, we require that $h$ behaves as $\mathcal{O}(k^{-3/2})$ whereas for piecewise cubic (P3) finite elements the criteria relaxes and $h$ should decrease as $\mathcal{O}(k^{-7/6})$. 
To find the best compromise between accuracy and computational efficiency, a common practice in finite-element simulation is to consider a fixed number of points $G$ per wavelength $\lambda = 2 \pi k^{-1}$ (i.e., $G \approx \lambda / h$), or simply taking $h$ decreasing as $\mathcal{O}(k^{-1})$. In practice, this is implemented  by using $h$-adaptive unstructured meshing.
In this study, we wish to use constant $G$ with $G=4$ to be consistent with the theoretical half-wavelength resolution of FWI \citep{Virieux_2009_OFW}. Indeed, the fill-in of the matrix increases as the polynomial degree increases. Therefore, a careful accuracy analysis of the FEFD method should be performed to find the smallest polynomial order providing the desired accuracy for $G$ = 4 as discussed in the {\it{Dispersion analysis}} section.

\subsubsection{Finite-difference frequency-domain (FDFD) discretization}
The second discretization scheme relies on the 27-point finite-difference stencil on regular Cartesian grid (constant grid interval $h$) in which compact $2^{nd}$-order accurate stencils minimize the numerical bandwidth and maximizes the sparsity of $A$ while reaching a high-order accuracy by mixing consistent mass and stiffness matrices on different (rotated) coordinate systems \citep{Operto_2007_FDFD,Brossier_2010_FNM,Chen_2012_DMF,Operto_2014_FAT,Aghamiry_2021_AFD}. These sparsity and compactness properties are useful to minimize the matrix fill-in induced by a sparse direct solver, which is used to solve the local problems in each subdomain of the preconditioner.
The stiffness and  consistent mass matrices are weighted by coefficients that are computed by least-squares minimization of the numerical dispersion and anisotropy. Classical implementation minimizes the dispersion jointly for a finite number of $G$, which means that the same coefficients are used in each row of the matrix  $A$ \citep{Operto_2007_FDFD,Brossier_2010_FNM,Chen_2012_DMF,Operto_2014_FAT}. Alternatively, $G$ can be first tabulated such that it covers the range of values found in the earth and the corresponding table of weighting coefficients are estimated once and for all by minimizing the phase-velocity numerical dispersion for each value of $G$ treated separately (the weights become functions of $G$). Then, for a simulation at a given frequency $f$, each row of the matrix $A$ (which is tied to a given grid point) is built by picking in the table the weights corresponding to the local $G=\lambda/h$ hence leading to the $\lambda$-adaptive 27-point stencil \citep{Aghamiry_2021_AFD}. We will use the latter in this study since \citet{Aghamiry_2021_AFD} showed that the $\lambda$-adaptivity increases the accuracy of the 27-point stencil quite significantly in heterogeneous media, in particular in presence of sharp contrast. Since the grid is regular, we will set $h$ such that $h=\lambda_{min}/4$ according to the FWI resolution criterion. 

In summary, the discretization rules used in FEFD and FDFD methods rely on two different paradigms. The discretization of the finite-element method is performed such that $G$ is roughly constant by adapting the size of the elements to the local wavelengths ($h$-adaptivity). By doing so, the accuracy of the numerical scheme is supposed to be similar all over the mesh. The discretization of the finite-difference method relies on a constant $h$ and hence $G$ varies from one grid point to the next. However, we locally adapt the weights of the stencil such that its accuracy is as uniform as possible as a function of $G$.

\subsubsection{Dispersion analysis of FEFD and FDFD methods}
\begin{figure}[ht!]
\begin{center}
\includegraphics[width=14cm,clip=true]{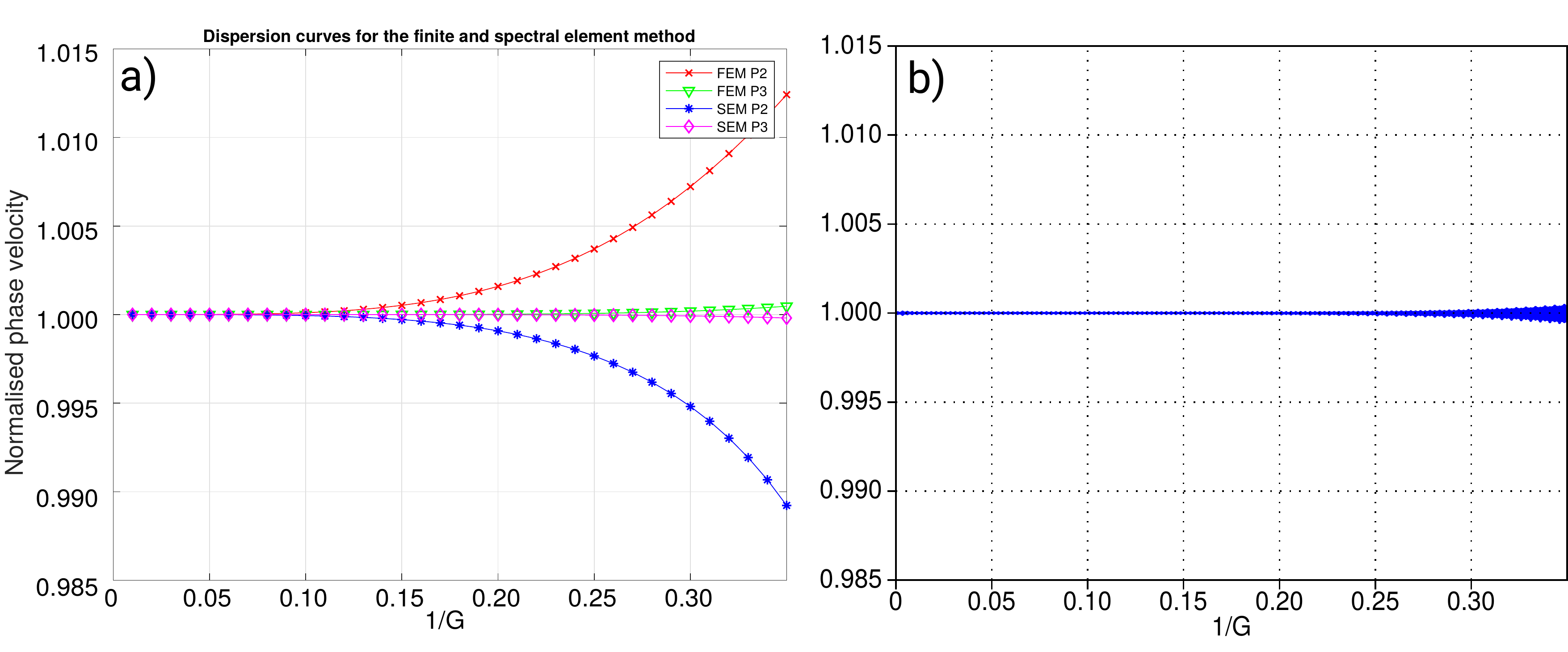}
\caption{Dispersion curves for (a) P2 and P3 spectral and classical finite elements, (b) the 27 FD scheme with adaptive coefficients.
The normalized phase velocity is the ratio between the numerical phase velocity and the wavespeed.}
\label{fig_disp}
\end{center}
\end{figure}

The results of the numerical dispersion analysis of the finite-element and finite-difference methods as developed by \citep{Ainsworth:2010:OBS} and \citet{Aghamiry_2021_AFD}, respectively, are shown in Fig.~\ref{fig_disp}. The figure shows the normalized numerical phase velocity against $1/G$ for plane waves in homogeneous media. The dispersion curve for the finite-element method suggests that the accuracy of P2 elements is insufficient for $G=4$, while a good accuracy is expected for P3 elements over a wide range of $G$ (Fig.~\ref{fig_disp}a). The dispersion curves of the $\lambda$-adaptive FDFD stencil suggests that its accuracy in homogeneous media is similar to that of the P3-FE stencil (Fig.~\ref{fig_disp}b).

\section{Numerical results}\label{sec:numerics}

\subsection{Benchmarks and experimental design}   
%
%
\subsubsection{Benchmarks}
We validate the FEFD/FDFD ORAS solver  against four benchmarks (Table~\ref{tab_bench_spec}). We limit ourselves to isotropic acoustic media with a constant density equal to 1. However, the conclusions of these tests should apply to visco acoustic media with heterogeneous density. All of the $\lambda$-adaptive FDFD simulations are performed with a discretization rule of four grid point per minimum wavelength, while the FEFD simulations are performed on $h$-adaptive tetrahedral mesh with four grid point per local wavelength, unless stated otherwise. The unstructured tetrahedral meshes are generated with the Mmg re-meshing software~\citep{dapogny2014three}. \\

{\it{Homogeneous and linear velocity models}} \\
The first two benchmarks involve an infinite homogeneous medium and an infinite laterally-homogeneous medium where the velocity linearly increases in the $y$ direction ($c(x,y,z)=c_0 + \alpha \times y$), referred to as gradient model. The size of both models is 10~km $\times$ 20~km $\times$ 10~km and the frequency is 7.5~Hz. The wavespeed in the  homogeneous medium is 1500~m/s. The grid interval in the Cartesian grid is 50~m, leading to $G$=4 for both FEFD and FDFD simulations (Table~\ref{tab_bench_spec}). In the gradient model, the wavespeed increases from 1500~m/s at y=0~km to 8500~m/s at y=20~km hence covering the wavespeeds encountered in the earth's crust and upper mantle (the velocity gradient is $0.35$ s$^{-1}$). The grid interval for the FDFD simulation is 50~m. Accordingly, $G$ increases from 4 to $\sim$23, while an $h$-adaptive mesh is used to perform FEFD simulation (Table~\ref{tab_bench_spec}). \\

The last two benchmarks involve the complex SEG/EAGE overthrust and GO\_3D\_OBS geomodels. \\

\begin{figure}[ht!]
\begin{center}
\includegraphics[width=16cm,clip=true]{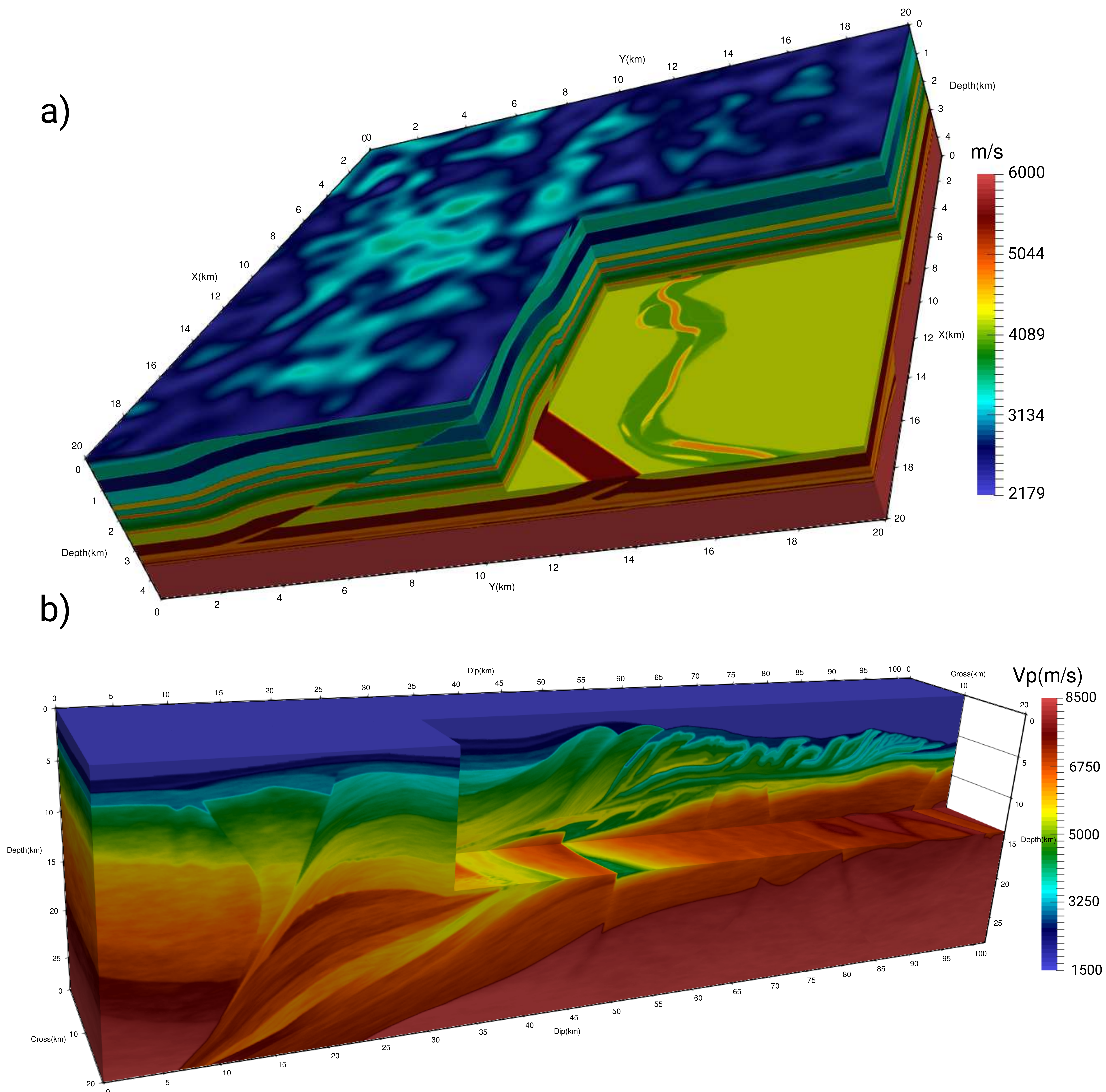}
\caption{The two realistic benchmarks used to assess the FEFD/FDFD ORAS-GMRES solver. (a) The 3D SEG/EAGE overthrust model \citep{Aminzadeh_1997_DSO}, representative of the exploration-geophysics scale. (b) Target of the regional GO\_3D\_OBS model representing the crust of a subduction zone \citep{Gorszczyk_2021_GNT}.}
\label{fig_over_go3dobs_model}
\end{center}
\end{figure}

\begin{figure}[ht!]
\begin{center}
\includegraphics[width=16cm]{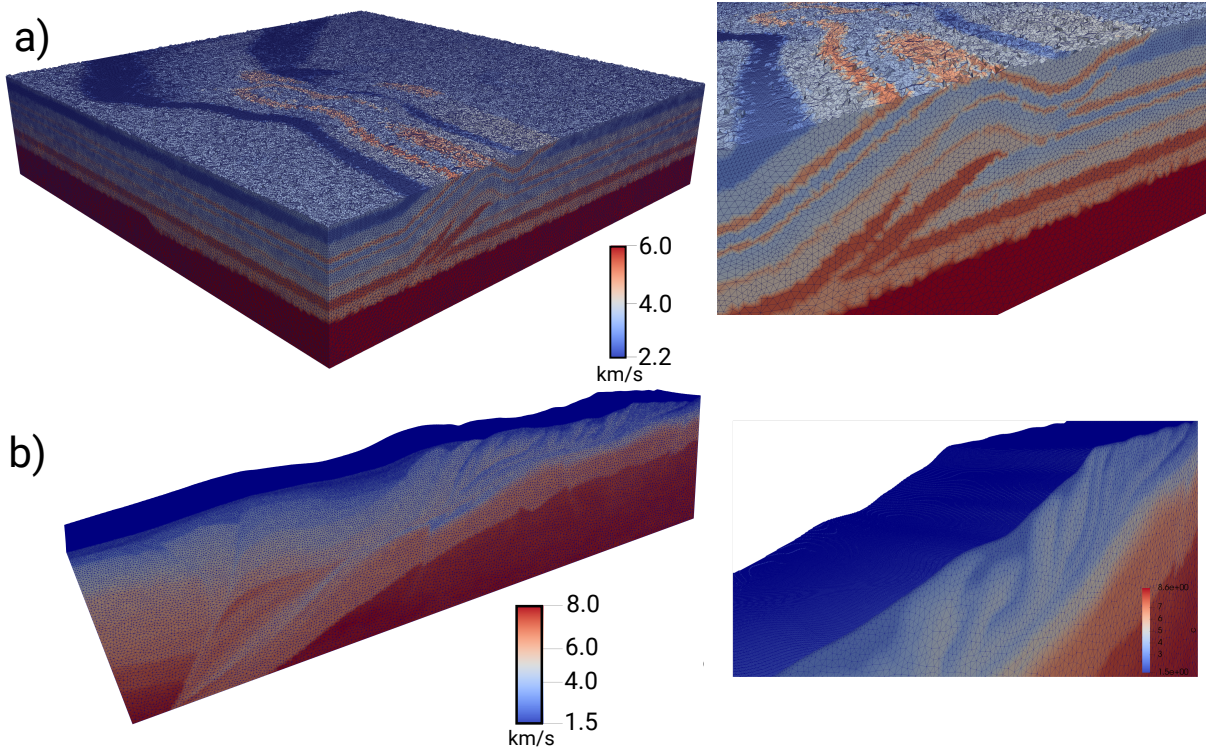}
\caption{Tetrahedral meshing of (a) the 3D SEG/EAGE overthrust and (b) GO\_3D\_OBS models. The right panels show a zoom.}
\label{fig_mesh_go3dobs}
\end{center}
\end{figure}

{\it{SEG/EAGE overthrust model}} \\
The overthrust model of dimensions 20~km $\times$ 20~km $\times$ 4.65~km is representative of the exploration geophysics scale, where $G$ typically ranges between 4 and 12 in a regular Cartesian grid. It represents a complex thrusted sedimentary succession on top of a structurally decoupled extensional and rift basement block \citep{Aminzadeh_1997_DSO} (Figure~\ref{fig_over_go3dobs_model}a). A complex weathered zone also characterizes it in the near-surface with sharp lateral velocity variations and several sand channels. The velocities range between 2179~m/s at the surface and 6000~m/s at the basement. The frequency is 10~Hz. For the FDFD simulation, we re-mesh the original model with a grid interval of 50~m and 8 grid points in the PMLs along each face of the 3D grid of dimensions 417 $\times$ 417 $\times$ 109 (Table~\ref{tab_bench_spec}). The unstructured tetrahedral mesh used for the FEFD simulation is illustrated in Figure~\ref{fig_mesh_go3dobs}(b-c). \\

{\it{GO\_3D\_OBS crustal geomodel}} \\
The 3D GO\_3D\_OBS crustal geomodel has been designed to assess seismic imaging techniques for deep crustal exploration  \citep{Gorszczyk_2021_GNT}. It covers a continental margin at a regional scale, where the wider range of wavespeeds found in the crust and upper mantle makes $G$ to range between 4 and 25 in a regular Cartesian grid. It embeds the main structural factors that characterize subduction zones and has been inspired by the FWI case study performed in the eastern Nankai trough by \citet{Gorszczyk_2017_TRW}. We select a target of the model of dimensions 20~km $\times$ 102~km $\times$ 28.4~km (Figure~\ref{fig_over_go3dobs_model}b). The simulations are performed for a frequency of 3.75~Hz except for the strong scaling analysis presented at the end of the section where a frequency as high as 10~Hz is processed. For the FDFD simulation, we discretize this target with a grid interval of 100~m with eight grid points in the PMLs, leading to a finite-difference grid of dimensions 217 $\times$ 1037 $\times$ 300, hence 67.5~million degrees of freedom. For the FEFD simulation, we build an unstructured mesh that complies with the bathymetry. Below the bathymetry, the size of the elements is set according to the local wavelength such that $G\approx4$ (Figure~\ref{fig_mesh_go3dobs}). 
\begin{table}
\begin{center}
\caption{Main specifications of the four benchmark models. $c_{m}$, $c_{M}$: Minimum and maximum wavespeeds. $f$: Frequency. $\lambda_{min}$,  $\lambda_{max}$: Minimum and maximum wavelength. $G_{min}$, $G_{max}$: Smallest and highest $G$. $N_\lambda$: Maximum number of propagated wavelengths.}
\label{tab_bench_spec}
\begin{tabular}{|c|c|c|c|c|c|c|c|c|c|}
\hline
Models    &        $c_{m}(m/s)$ 	&  $c_{M}(m/s)$	& $f(Hz)$  & $\lambda_{m}(m)$ & $\lambda_{M}(m)$ & $G_{m}$ & $G_{M}$ & $N_\lambda$		\\ \hline
Homogeneous &  1500 & 1500 & 7.5  & 375	& 375 & 4 &  4	&  40	\\ \hline
Linear 		 &  1500 & 8500 & 7.5 &  375	& 1133 & 4 & 22.7 &	29	\\ \hline
Overthrust 	 & 2179 & 6000 & 10 &  	218 & 600 & 4.4 & 12 &	50		\\ \hline
GO\_3D\_OBS 	 & 1500 & 8639.1 & 3.75 &  400	& 2303.8 & 4 & 23 &	255		\\ \hline
\end{tabular}
\end{center}
\end{table}
%
\subsubsection{Experimental setup}
For the first two benchmarks, the accuracy of the FEFD and FDFD solutions are assessed against analytical solutions. The reader is referred to \citet{Kuvshinov_2006_EST} for the analytical solution of the Helmholtz equation in a velocity gradient model. \\
For the last two benchmarks, the FEFD and FDFD solutions are assessed against a reference solution computed with the highly-accurate Convergent Born Series (CBS) method \citep{Osnabrugge_2016_CBS}. Compared to finite-difference or finite-element methods that rely on discretization of differential operators, the CBS method belongs to another class of methods referred to as volume integral methods based upon the Green's function theorem. More specifically, the solution of the Helmholtz equation in heterogeneous media is inferred from Fourier-domain analytic expression of the Green's functions in homogeneous media by solving the scattered-field wave equation with material contrasts of arbitrary scattering strength starting from the zero-order term of the Born series and iteratively recovering the higher-order scattering terms.
Therefore, we expect to achieve an accuracy of the CBS solutions at the level of machine precision while \citet{Osnabrugge_2016_CBS} conclude that the CBS method is nine orders more accurate than the pseudo-spectral time-domain method. For all the extracted wavefields using CBS, the stopping criterion of iteration is set according to a backward error of $10^{-12}$. \\
In addition to visual comparison, we assess the accuracy of the FEFD and FDFD stencils by computing the $\ell{1}$ norm of the difference between the CBS and FEFD/FDFD wavefields (field Err in Table~\ref{tab_bench_results}) according to the formula
\begin{equation}
\text{Err}=\frac{| \bold{W}  \mathcal{R}(\bold{u}_{cbs}-\bold{u}_{fe/fd}) |}{ | \bold{W}  \mathcal{R}(\bold{u}_{cbs})|}+\frac{|  \bold{W} \mathcal{I}(\bold{u}_{cbs}- \bold{u}_{fe/fd}) |}{  | \bold{W} \mathcal{I}(\bold{u}_{cbs}) |},
\label{eqerror}
\end{equation}
where  $\bold{u}_{cbs}$ and $\bold{u}_{fe/fd}$ denotes the CBS and FEFD/FDFD wavefields, respectively, $\mathcal{R}$ and $\mathcal{I}$ denote the real and imaginary parts of a complex number, respectively, and $\bold{W}$ is a diagonal matrix applying a linear gain with distance from the source to the wavefield for amplitude balancing. We also mute the contribution of the nearest points to the source to avoid introducing some bias generated by the different implementation of a point source in the FEFD and FDFD schemes. \\
%
%
The source is a point source located on a grid point where the temporal source signature is a delta function: $s(\bold{x},\omega)= \delta(\bold{x}-\bold{x}_s)$, where $\delta$ denotes the delta function and $\bold{x}_s$ denotes the source coordinates. We implement PMLs along each face of the grid. \\
We perform the simulation on the Occigen supercomputer of CINES (\url{https://www.cines.fr}). The computer nodes contain two  Haswell E5-2690V3@2.6 GHz processors with 128 Giga bytes of shared memory and 12 cores per processor. The high-speed network is Infiniband FDR 56 Gbit/s.
\subsection{Tuning of the ORAS solver}
Before showing the results of the ORAS solver for the four benchmarks when the subsurface model is discretized with the FEFD and FDFD methods, we tune key ingredients of the ORAS solver such that the best compromise between accuracy and computational efficiency is reached. These ingredients include the stopping criterion of iterations, arithmetic precision (single versus double precision), the solver used to factorize local matrices and the orthogonalization scheme of the Krylov basis. We perform this tuning with the 3D GO\_3D\_OBS benchmark and the  FDFD scheme. 

\subsubsection{Stopping criterion of GMRES iterations} 
\begin{figure}[ht!]
\begin{center}
\includegraphics[width=16cm,clip=true]{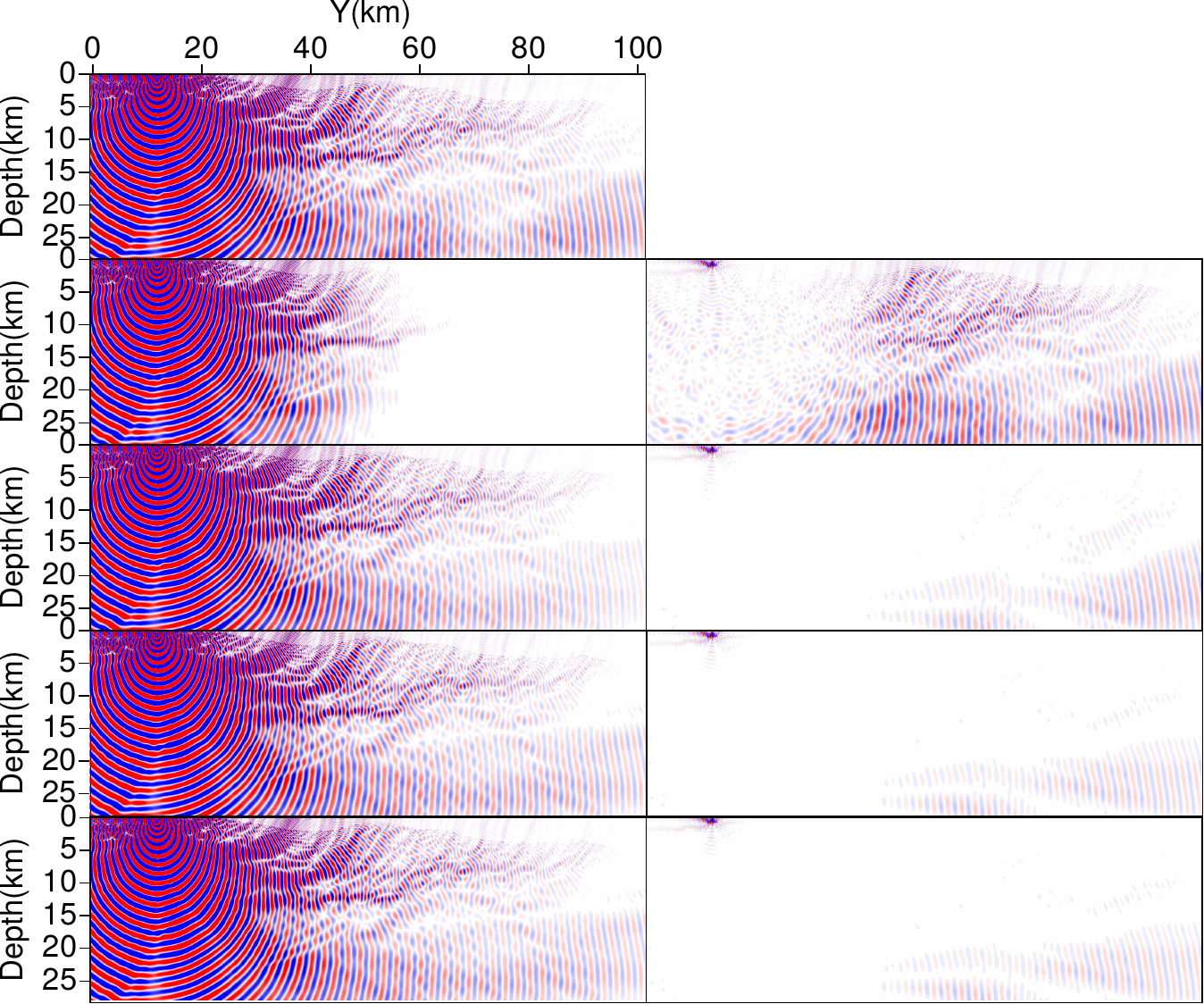}
\caption{Setting the stopping criterion of GMRES iterations. Frequency is 3.75~Hz. In the caption of this figure, the panels are labeled with matrix notations: $a_{i,j}$ where $i$ and $j$ denote the row and column indexes, respectively. $a_{1,1}$: Section of the CBS wavefield (real part). $a_{2:5,1}$: From top to bottom, FDFD wavefields (real part) for $\varepsilon_{gmres}= 10^{-2}, 10^{-3}, 10^{-4}, 10^{-5}$. $a_{2:5,2}$: Differences between CBS and FDFD wavefields.}
\label{fig_wavefield_tolerance}
\end{center}
\end{figure}

\begin{figure}[ht!]
\begin{center}
\includegraphics[width=13cm,clip=true]{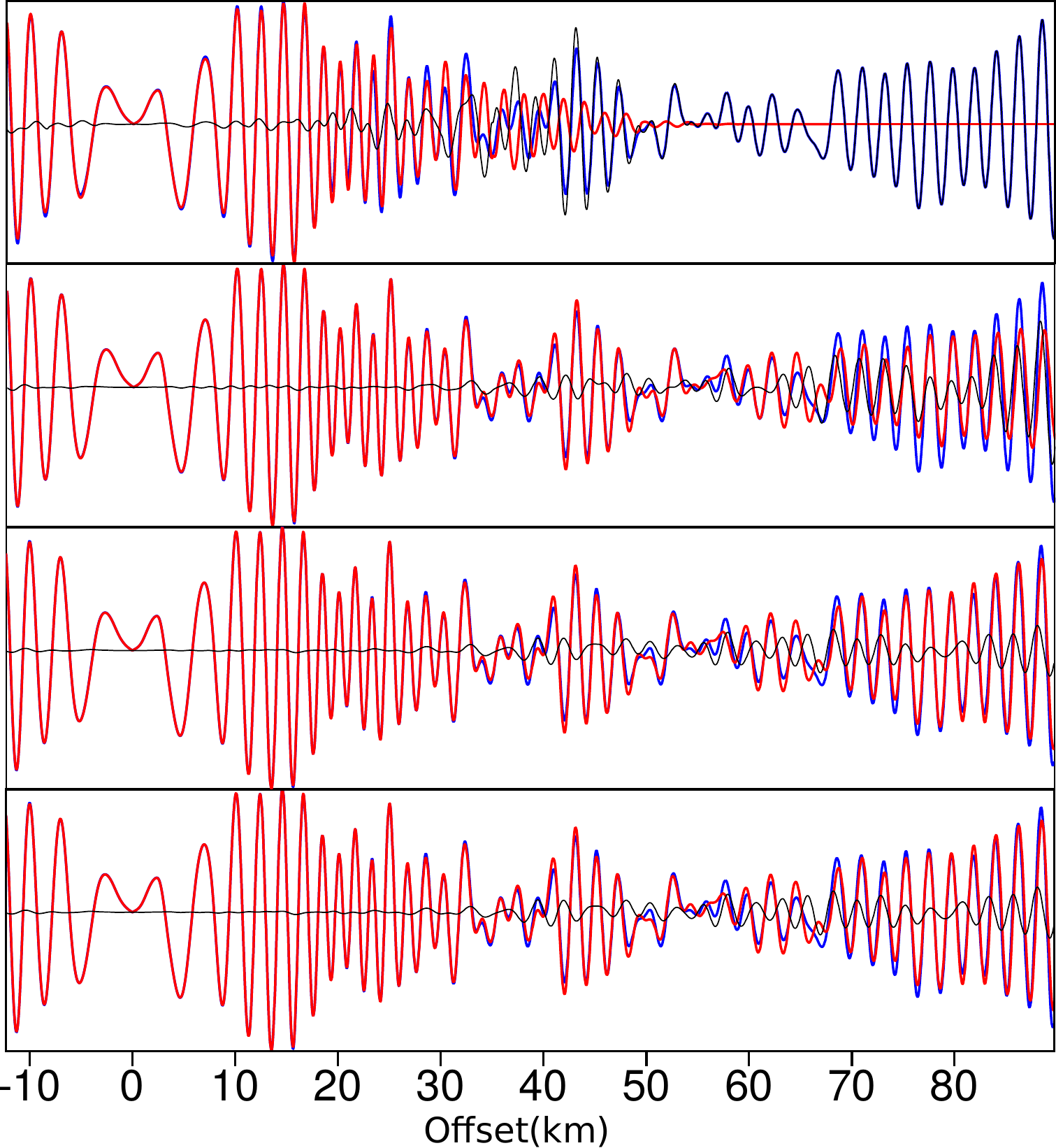}
\caption{Stopping criterion of GMRES iterations. From top to bottom, direct comparison along a y-profile at 16~km depth (x=10.5~km) between CBS wavefield (blue) and FDFD wavefields (red) computed with $\varepsilon_{gmres}=  10^{-2}, 10^{-3},  10^{-4}$, $10^{-5}$ (Figure~\ref{fig_wavefield_tolerance}). The differences are plotted in black. A linear gain with offset is applied for amplitude balancing.}
\label{fig_log_tolerance}
\end{center}
\end{figure}

The simulations are performed in single precision on 660~cores.
As stopping criterion of iteration, we use the relative backward error defined by
\begin{equation}
\varepsilon_{gmres}=\frac{\|A u - f \|_2^2}{\|f\|_2^2}.
\end{equation}
We test $\varepsilon_{gmres}=10^{-2}$, $10^{-3}$, $10^{-4}$ and $10^{-5}$ (Table~\ref{tab_epsilon} and Figures~\ref{fig_wavefield_tolerance} and \ref{fig_log_tolerance}). The number of iterations required to satisfy this criterion is 13, 32, 45 and 61, respectively.
The errors of the FDFD simulations performed with $\varepsilon_{gmres}=10^{-4}$ and $\varepsilon_{gmres}=10^{-5}$ are similar and hence the former value is used in the remaining of the paper. Using $\varepsilon_{gmres}=10^{-4}$  instead of $\varepsilon_{gmres}=10^{-5}$ reduces the number of iterations from 61 to 45 and the computational time by 27$\%$.

\begin{table}[ht!]
\begin{center}
\caption{Cost of a simulation in the GO\_3D\_OBS model as a function of the GMRES backward error criterion. $\#it$: number of iterations. $T_{s}(s)$: Time elapsed in GMRES to compute the solution. $Err$: Error between the CBS and FDFD wavefields as defined by equation~\ref{eqerror}.}
\label{tab_epsilon}
\begin{tabular}{|c|c|c|c|}
\hline
\multicolumn{4}{|c|}{\bf{Stopping criterion of GMRES iteration}} \\ \hline
$\varepsilon_{gmres}$  & $\#it$ & $T_{s}(s)$ & Err \\ \hline
$10^{-2}$     &    13  &  3.5 & 0.8663 \\ \hline
$10^{-3}$     &    32  &  8.2 & 0.2567 \\ \hline
$10^{-4}$     &    45  & 11.5 & 0.2096 \\ \hline
$10^{-5}$     &    61  & 15.9 & 0.2093 \\ \hline
\end{tabular}
\end{center}
\end{table}

\subsubsection{Arithmetic precision}
The results shown in the section {\it{FEFD versus FDFD discretization: accuracy and computational cost}} are performed in single precision. We check here that a similar accuracy is reached when we perform simulation in single and double precisions. Compared to the previous section, the simulations are performed on 1344 cores because the double precision simulation runs out of memory on 660 cores. Single precision allows for the solution time to be decreased by a factor two compared to double precision simulations (fourth column of Table~\ref{tab_arithmetic}) while the error between the FDFD and the CBS simulations is the same, hence validating the single precision option (last column of Table~\ref{tab_arithmetic}).
\begin{table}[ht!]
\begin{center}
\caption{Footprint of arithmetic precision on computational cost and accuracy. \#it: Number of iterations. $T_{f}(s)$: Elapsed time for local LU factorizations. $T_{s}(s)$: Elapsed time for GMRES. $Err$: Error between the CBS and FDFD wavefields as defined by equation~\ref{eqerror}.}
\label{tab_arithmetic}
\begin{tabular}{|c|c|c|c|c|}
\hline
\multicolumn{5}{|c|}{\bf{Arithmetic precision}} \\ \hline
Precision  & $\#it$ & $T_{f}(s)$ & $T_s(s)$ & Err \\ \hline
Double     & 53 & 7.9 & 14.0 &  0.2097 \\ \hline
Single     &   53 &  4.4 &  7.1  & 0.2096  \\ \hline
\end{tabular}
\end{center}
\end{table}

\subsubsection{Solution of the local problems and multiple right-hand sides}
In this section, we assess two direct solvers and orthogonalization schemes of the Krylov basis.
We test the Intel MKL PARDISO solver \citep{Bollofer_2019_LSI,Bollofer_2020_SSD,Alappat_2020_RAC} and the MUMPS solver \citep{Amestoy_2018_ESM,Amestoy_2019_EUS,MUMPS_2021_MMP} to factorize the local matrices tied to each subdomain. When  using MUMPS, we test both the full-rank (FR) and the block low-rank (BLR) version of the solver, referred to as MUMPS$_{FR}$ and MUMPS$_{BLR}$, respectively \citep{Amestoy_2016_FFF,Amestoy_2019_PSB}.
In the latter case, we use as compression threshold $\varepsilon_{BLR}=10^{-3}$. We also compare the performances when a classical Gram-Schmidt (CGS) and a modified Gram-Schmidt (MGS) algorithms are used for the orthogonalization of the Krylov basis \citep{Saad_2003_IMS} and we perform the computations in single and double precisions. Moreover, we perform the simulation with one RHS and 130 RHSs, the latter mimicking a sparse seabed node acquisition. For multi-RHS simulation, we test a pseudo-block Krylov method, a block Krylov method and block Krylov with recycling by processing the 130 RHSs in group of 20 \citep{Parks_2006_RKS,Jolivet_2016_BIM}.
The results are outlined in Table~\ref{tab_solver}. The smaller computational time is reached with MUMPS$_{BLR}$ although the BLR approximation generates some iteration overheads. The CGS orthogonalization provides substantial computational saving compared to the MGS counterpart in case of multi-RHS simulation. These simulations further confirm that single precision arithmetic decreases by a factor almost two the cost of the simulation compared to double precision counterpart. Concerning multi-RHS simulation, the best performance is obtained with the pseudo-block Krylov method with a speedup of $7.2 \times 130 / 256.8 = 3.64$ compared to the case where each RHS is processed separately (last row of Table~\ref{tab_solver}). This suggests that the most basic pseudo-block Krylov remains the most efficient one for seismic imaging applications, which involve moderate material contrast and simple wave propagation phenomena. We check that for groups of RHS larger than 20 no improvement is achieved during multi-RHS simulation.
\begin{table}[ht!]
\begin{center}
\caption{$Solver$: Intel MKL Pardiso versus MUMPS$_{FR/BLR}$. $arith$: $single$ versus $double$ precision. $ortho$: CGS versus MGS orthogonalization. \#it: Number of iterations. $T_{f}(s)$: Elapsed time for local LU factorizations. $T_{s}(s)$: Elapsed time for all GMRES iterations for one RHS and 130 RHSs processed with a pseudo-block Krylov method. $T_{tot}(s)=T_{f}(s)+T_{s}(s)$: Total elapsed time for the simulation.}
\label{tab_solver}
\begin{tabular}{|c|c|c|c|c|c|c|c|c|c|}
\hline
& & & & \multicolumn{3}{|c|}{1 RHS} & \multicolumn{3}{|c|}{130  RHS} \\ \hline
Solver   &  arith. &  ortho. & $T_{f}(s)$ & \#it & $T_s(s)$ & $T_{tot}(s)$ & \#it & $T_s(s)$ & $T_{tot}(s)$ \\ \hline
PARDISO  & double & MGS & 12.1 & 48 & 16.5 & 28.6 & ~50 & 724.4 & 736.5 \\ \hline
PARDISO  & double & CGS & 12.1 & 48 & 16.7 & 28.8 & ~50 & 575.5 & 587.6 \\ \hline
MUMPS$_{FR}$    & double & MGS & 11.3 & 48 & 17.2 & 28.5 & ~50 & 654.9 & 666.2 \\ \hline
MUMPS$_{FR}$    & double & CGS & 11.3 & 48 & 16.9 & 28.2 & ~50 & 482.2 & 493.5 \\ \hline
PARDISO  & single & MGS &  6.5 & 48 &  9.9 & 16.4 & ~50 & 403.6 & 410.1 \\ \hline
PARDISO  & single & MGS &  6.5 & 48 &  9.6 & 16.1 & ~50 & 334.8 & 341.3 \\ \hline
MUMPS$_{FR}$    & single & MGS &  6.2 & 48 &  8.9 & 15.1 & ~50 & 345.8 & 352.0 \\ \hline
MUMPS$_{FR}$    & single & CGS &  6.2 & 48 &  8.9 & 15.1 & ~50 & 275.3 & 281.5 \\ \hline
MUMPS$_{BLR}$ & single & CGS &  4.9 & 53 &  7.2 & 12.1 & ~55 & 256.8 & 261.7 \\ \hline
\end{tabular}
\end{center}
\end{table}
%
%

\subsection{FEFD versus FDFD discretization: accuracy and computational cost}
According to the results of the previous section, we assess the FDFD and FDFD schemes with $\varepsilon_{gmres}=10^{-4}$. Computations are performed in single precision. The local solver relies on MUMPS$_{BLR}$ with $\varepsilon_{BLR}=10^{-3}$ while a CGS orthogonalization is used to build the Krylov basis. For FDFD simulation, we use the one-level ORAS preconditioner, while the two-level counterpart is used for the FEFD simulations. PML conditions are used between the subdomains in the FDFD simulation with an overlap of 3, while a Robin condition is used for the FEFD simulation. The results of the numerical tests for the four benchmarks are outlined in Table~\ref{tab_bench_results}.
\begin{table}
\begin{center}
\caption{Results of the four benchmarks obtained with the FEFD and FDFD methods. $\#d$: Number of dofs (including PMLs). $\#cores$: Number of cores. $\#it$: Number of GMRES iterations (for FEFD, average number of inner iterations to solve the coarse problem in parentheses). $T_{f}(s)$: Elapsed time for local LU factorizations. $T_{s}(s)$: Elapsed time for all GMRES iterations. $T_{tot}(s)=T_{f}(s)+T_{s}(s)$: Total elapsed time for the simulation.  $T_{hc}(s) = \#c \times T_{tot}$: Scalar time. ERR: Solution error as defined by equation~\ref{eqerror}. The specifications of the four benchmarks in term of size and frequency are outlined in Table~\ref{tab_bench_results}.}
\label{tab_bench_results}
\begin{tabular}{|c|c|c|c|c|c|c|c|c|c|c|c|c|c|c|}
\hline
              &   \multicolumn{8}{|c|}{FEFD method}                                                               \\ \hline
Benchmark       &    $\#d(M)$            & \#cores       &  \#it        & $T_{f}(s)$     & $T_{s}(s)$  & $T_{tot}(s)$ & $T_{hc}(h)$ &  Err      		\\ \hline
Homogeneous     &      526.2          &  2400     &  31(13)            &   64.1      &    73.0     &  137.1 & 91.4 & 0.0891        \\ \hline
Gradient 	                 &     147.1          &  2400        &	5(14)       &  12.4       &    6.1   &  18.5    &  12.3 & 0.0485       \\ \hline
Gradient$^{(c)}$ 	          &    526.2          &  2400        &	8(16)       &    59.8     &    22.3  & 82.1      & 54.7  & 0.0119    \\ \hline
Overthrust                   &       157.1          &  2400        &	    6(17)   &    13.5     &   6.7   &   20.2   & 13.5 & 0.5215            \\ \hline
Overthrust$^{(c)}$             &     516.5          &  2400        &	    6(20)   &  57.6       &     19.5  &   77.1   & 51.4 & 0.5156             \\ \hline
GO\_3D\_OBS                      &    597.7         &  2400        &	10(14)	    &    62.3     &   26.1  &     88.4    &  58.9 & 0.4212                     \\ \hline
\hline
                            &             \multicolumn{8}{|c|}{FDFD method}                                                \\ \hline
Benchmark            & $\#d(M)$            & \#c       &  \#it        & $T_{f}(s)$     & $T_{s}(s)$ & $T_{tot}(s)$ & $T_{hc}(h)$  & Err      		\\ \hline
Homogeneous              &       19.6           & 396   &  41     & 5.4     & 5.7  &  11.1   & 1.2 & 0.0317     \\ \hline
Gradient 	           &     23.5           & 396   &  42     & 5.7     & 6.4  &  12.1   & 1.3 & 0.0114     \\ \hline
Overthrust             &       19.1           & 363   &  33     & 5.0      &  4.5    & 9.5 &  0.96  & 0.1188            \\ \hline
GO\_3D\_OBS               &        67.5           & 660   &  45     & 10.8     &   11.5   & 22.3   & 4.1 & 0.2096                     \\ \hline
\end{tabular}
\end{center}
\end{table}

\subsubsection{Homogeneous medium}
\begin{figure}[htb!]
\begin{center}
\includegraphics[width=16cm,clip=true]{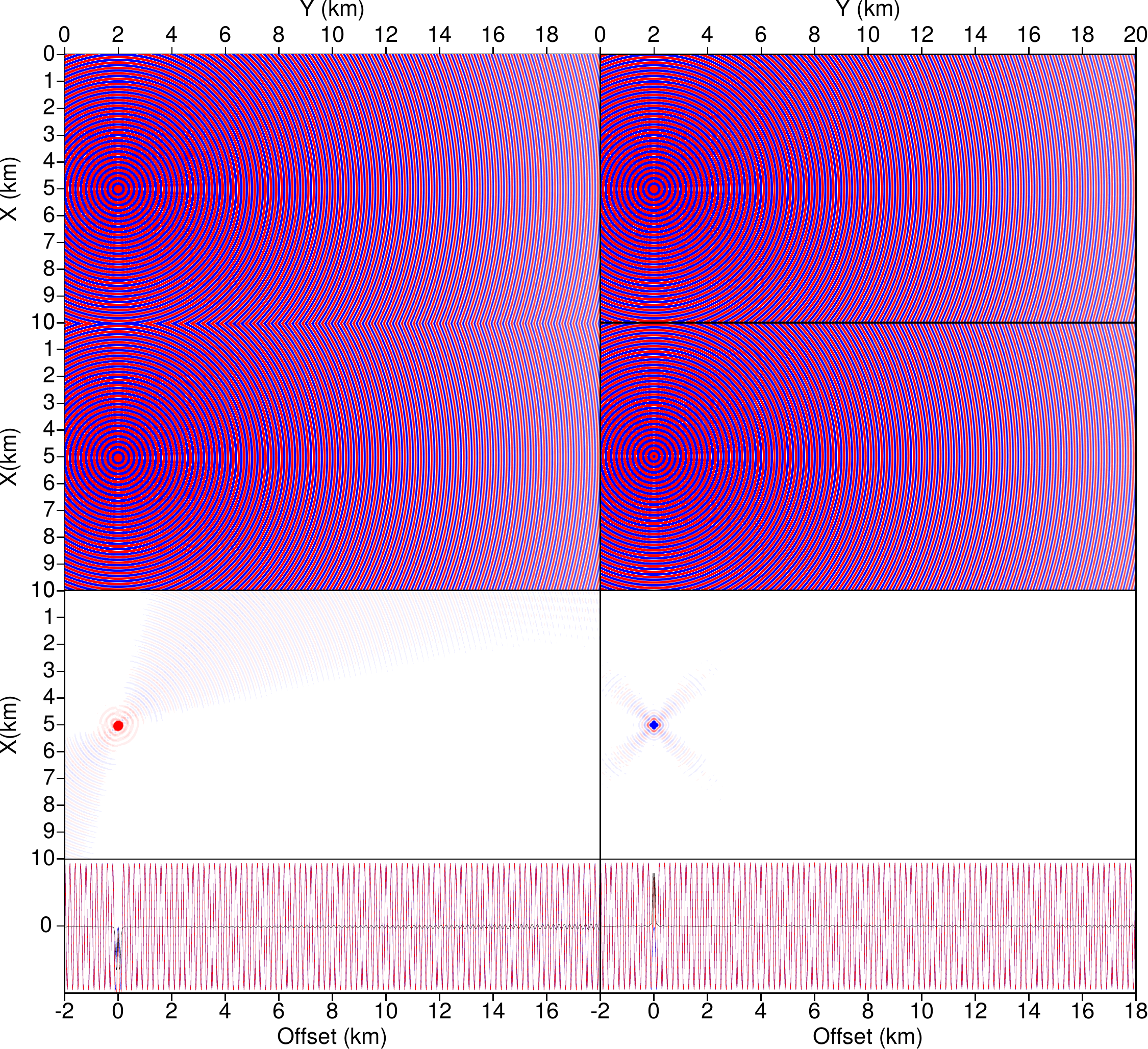}
\caption{FE and FD monochromatic wavefields in infinite homogeneous medium. $a_{1,1:2}$: Analytical solution (real part). $a_{2,1}$: FEFD wavefield (real part). $a_{2,2}$: FDFD solution (real part). $a_{3,1}$ : Difference between Analytical solution and FEFD wavefield. $a_{3,2}$:  Same as $a_{3,1}$ for the FDFD wavefield.  $a_{4,1}$: Direct comparison between analytical (blue) and FEFD solution (red) along a $Y$ profile cross-cutting the source position. Difference is plot in black. $a_{4,2}$: Same as $a_{4,1}$ for the FDFD wavefield.}
\label{fig_homogeneous_wavefield_log}
\end{center}
\end{figure}

The source is positioned at (5~km, 2~km, 5~km). We perform the FEFD and FDFD simulations with 2400 and 396 cores, respectively. The scalar time (the parallel run time $\times$ the number of cores, where the parallel run time is the sum of the elapsed time required to perform LU factorization in the subdomains and the GMRES iterations) is 91.4 hours and 1.2 hour for the FEFD and FDFD simulations, respectively (Table~\ref{tab_bench_results}). The analytical and the numerical FEFD and FDFD wavefields in the (x,y) plane across the source and the difference between the analytical and the numerical solutions are shown in Figure~\ref{fig_homogeneous_wavefield_log} as well as the direct comparison between the analytical and numerical wavefields along the $y$-profile. The results show a slightly better accuracy of the FDFD solution relative to the FEFD counterpart (Table~\ref{tab_bench_results}). However, the two solutions are quite close, hence supporting the conclusions of the dispersion analysis (Fig~\ref{fig_disp}).
%
%

%
%
\subsubsection{Velocity gradient model}
\begin{figure}[htb!]
\begin{center}
\includegraphics[width=16cm,clip=true]{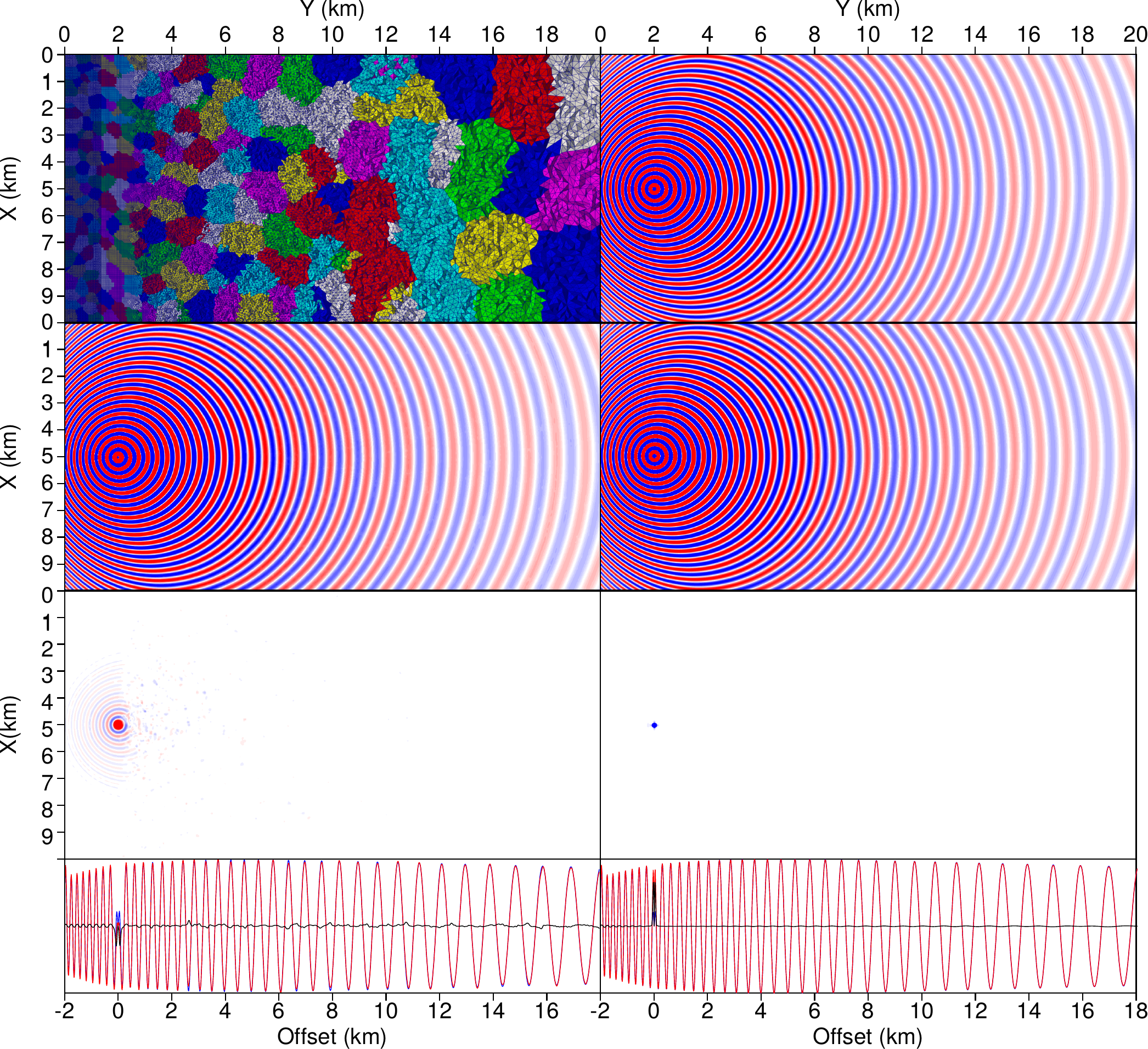}
\caption{$h$-adaptive FE and $\lambda$-adaptive FD monochromatic wavefields in linear velocity model. $a_{1,1}$: $h$-adaptive tetrahedral mesh for FDFD simulation. Each subdomain is plotted with a different color. $a_{1,2}$: Analytical solution (real part). $a_{2,1}$: FEFD wavefield (real part). $a_{2,2}$: FDFD wavefield (real part). $a_{3,1}$: Difference between analytical and FEFD solution. $a_{3,2}$:  Same as $a_{3,1}$ for the FDFD solution. $a_{4,1}$: Direct comparison between analytical (blue) and FEFD (red) solution along a $Y$ profile cross-cutting the source position. Difference is plot in black. $a_{4,2}$: Same as $a_{4,1}$ for FDFD solution.}
\label{fig_gradient_wavefield_log}
\end{center}
\end{figure}

\begin{figure}[htb!]
\begin{center}
\includegraphics[width=16cm,clip=true]{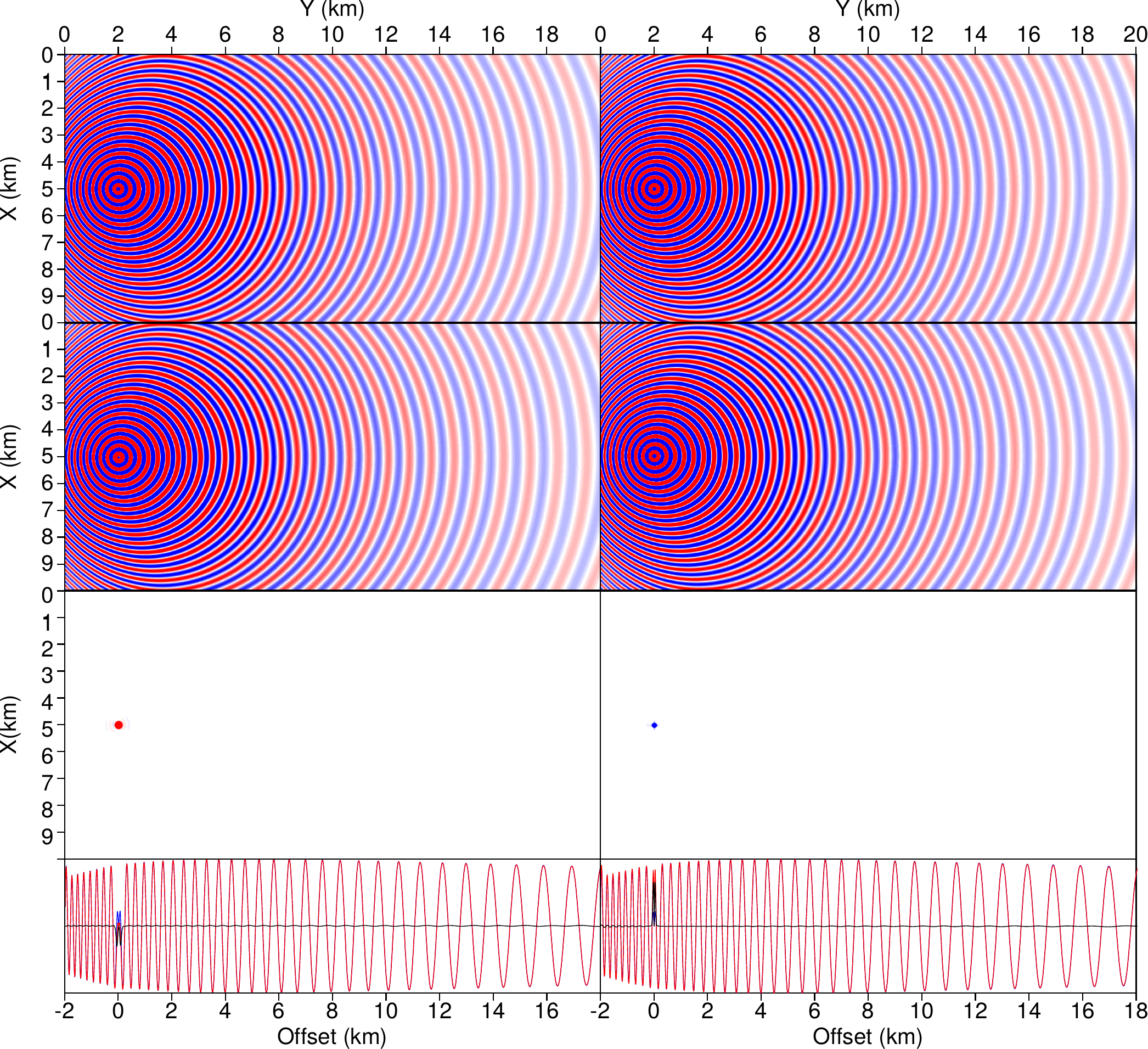}
\caption{Same as Figure~\ref{fig_gradient_wavefield_log} when FEFD is performed on a regular mesh with a discretization rule of four grid point per minimum wavelength.}
\label{fig_gradient_wavefield_log_car}
\end{center}
\end{figure}

Our motivation behind this test, beyond accuracy issue, is to gain a first insight whether the $h$-adaptivity implemented in the FEFD method can balance the highest number of dofs per element relative to the FDFD method in term of computational cost when a wide range of wavelengths are involved in the simulation.  
The results are shown in Figure~\ref{fig_gradient_wavefield_log} with the same showing as for the homogeneous case. As for the homogeneous case, the FDFD method remains marginally more accurate  and remains one order of magnitude less expensive than the FEFD counterpart (Check field $T_{hc}$ in Table~\ref{tab_bench_results}). 
Figure~\ref{fig_gradient_wavefield_log_car} shows the results of the FEFD simulation in a regular mesh where the number of grid point per minimum wavelength is set to four. The accuracy is slightly improved at the expense of the computational cost, which increases by a factor 4.4 (Table~\ref{tab_bench_results}).
\subsubsection{3D SEG/EAGE overthrust velocity model}
\begin{figure}[htb!]
\begin{center}
\includegraphics[width=16cm,clip=true]{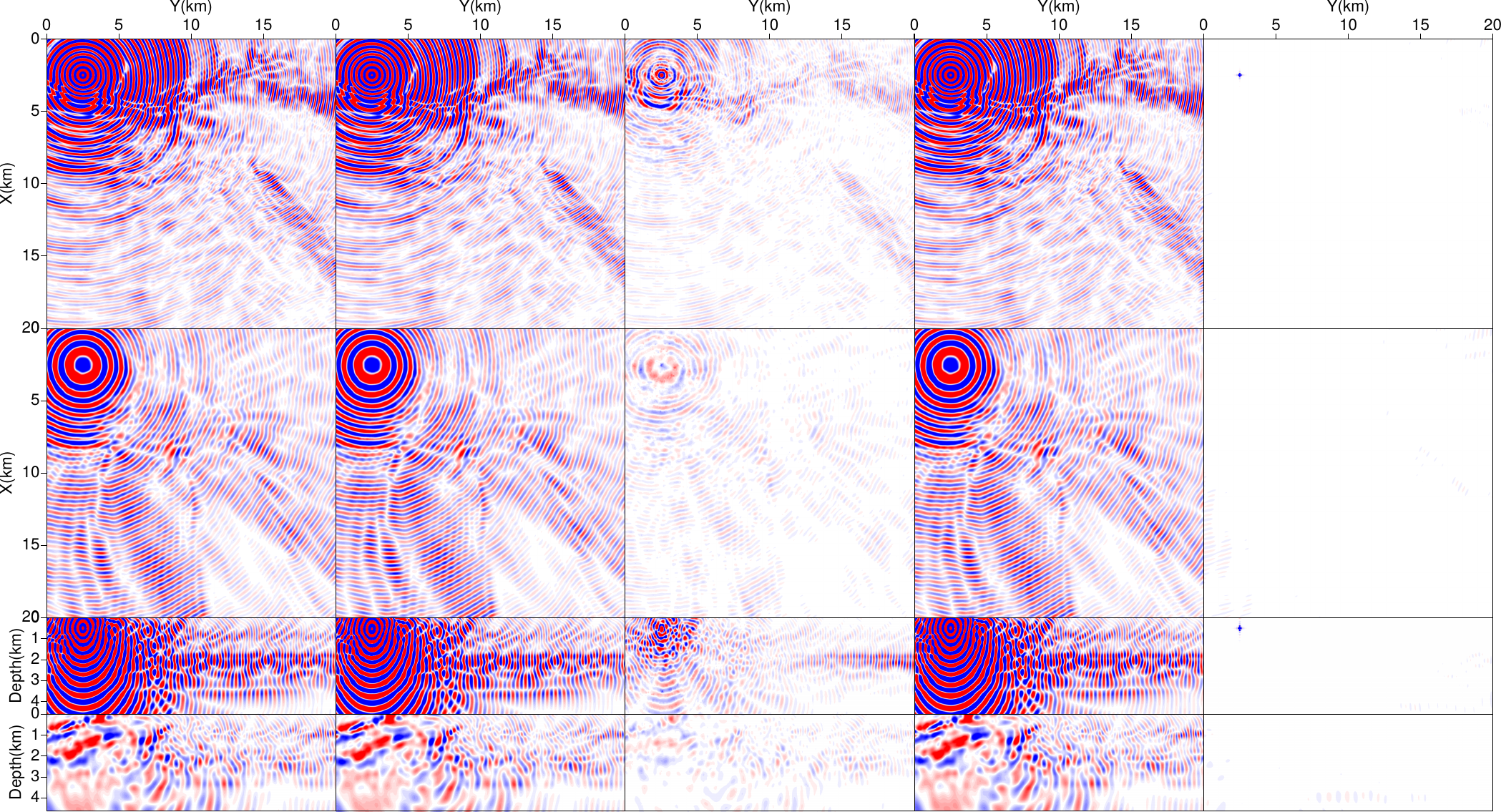}
\caption{3D EAGE/SEG overthrust model. Comparison between the wavefields computed with the CBS, FEFD and FDFD methods. The rows shown from top to bottom two depth slices at 500~m depth (across the source) and 2~km depth, and two vertical sections at x=2.5~km (across the source) and 15~km. From left to right, the columns show the CBS wavefield, the FEFD wavefield, the differences between the two, the FDFD wavefield and its differences with the CBS wavefield.}
\label{fig_over_wavefield}
\end{center}
\end{figure}

\begin{figure}[htb!]
\begin{center}
\includegraphics[width=16cm,clip=true]{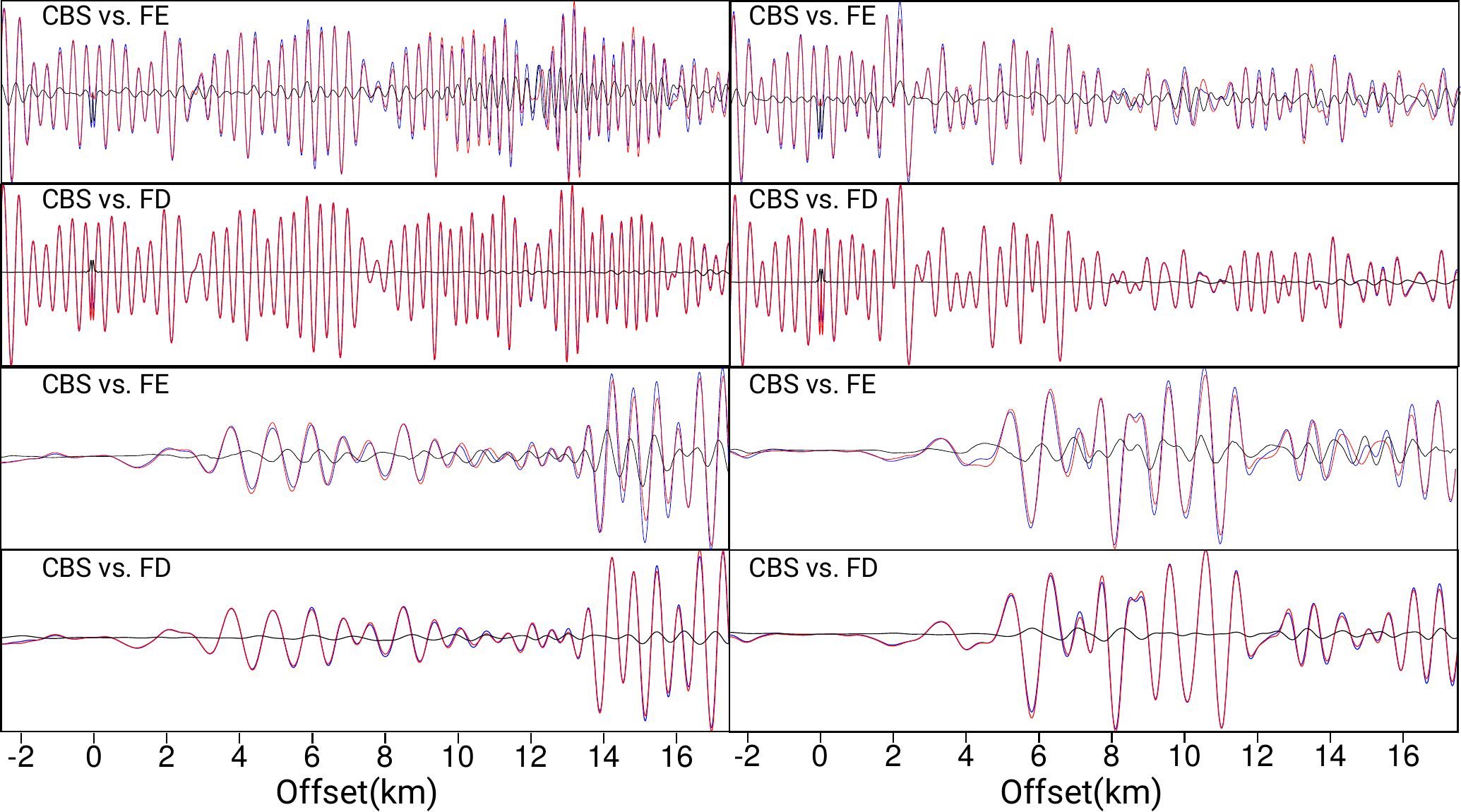}
\caption{3D EAGE/SEG overthrust model. Direct comparison between the CBS (blue) and FEFD/FDFD wavefields (red). The differences are plotted in black. ($a_{1:2,1}$): Profile at (x,z)=(2000~m, 500~m). ($a_{1:2,2}$): Profile at (y,z)=(2500~m, 500~m). ($a_{3:4,1}$): Profiles at (x,z)=(15000~m, 3500~m). ($a_{3:4,2}$): Profile at (y,z)=(15000~m, 3500~m). ($a_{1,1:2}, a_{3,1:2}$): FEFD solution. ($a_{2,1:2}, a_{4,1:2}$): FDFD solution.}
\label{fig_over_log}
\end{center}
\end{figure}

\begin{figure}[htb!]
\begin{center}
\includegraphics[width=16cm,clip=true]{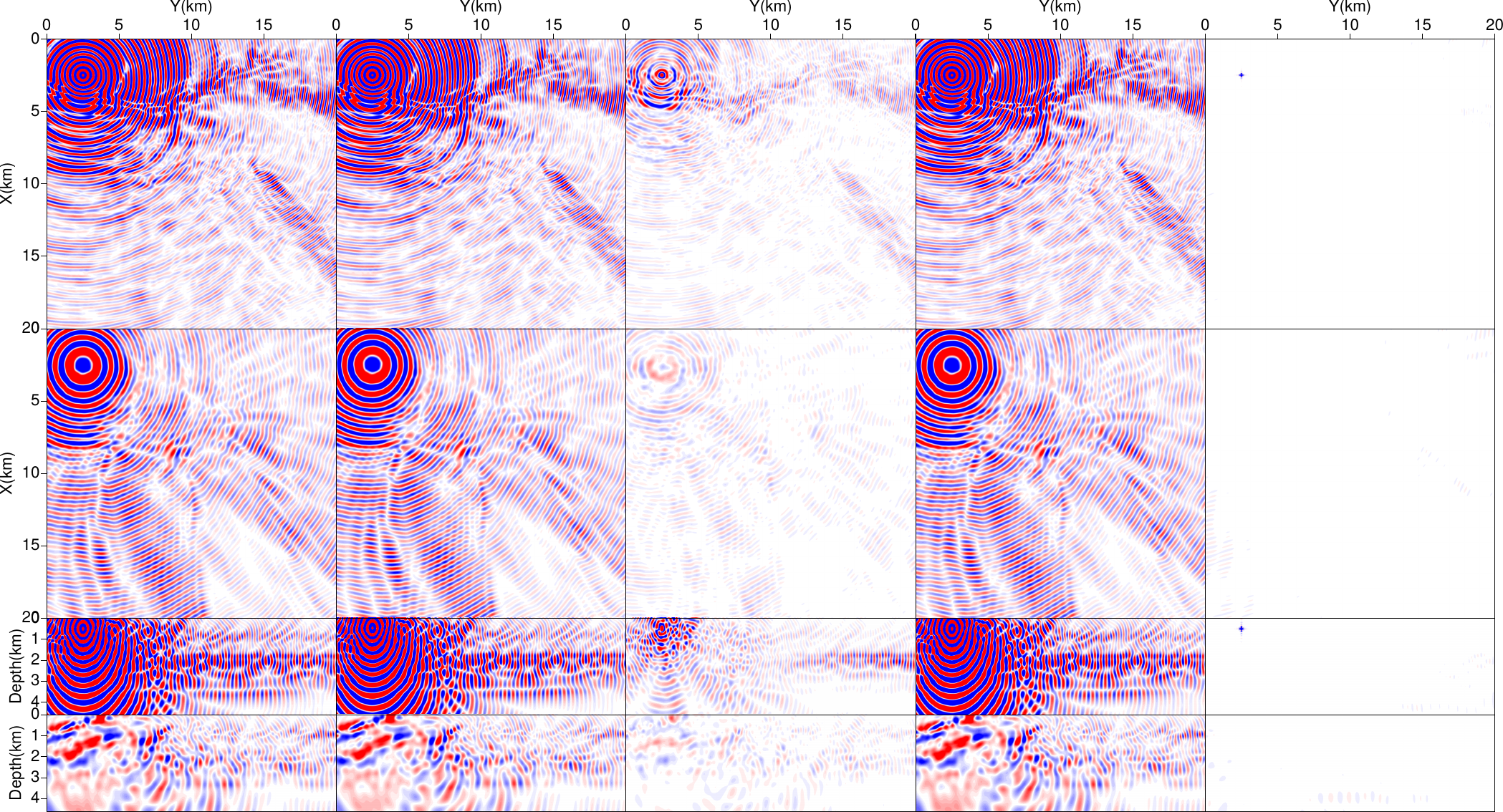}
\caption{3D EAGE/SEG overthrust model. Same as Figure~\ref{fig_over_wavefield} but the FEFD simulation is performed on regular mesh with a discretization rule of four grid point per minimum wavelength.}
\label{fig_over_wavefield_car}
\end{center}
\end{figure}

\begin{figure}[htb!]
\begin{center}
\includegraphics[width=16cm,clip=true]{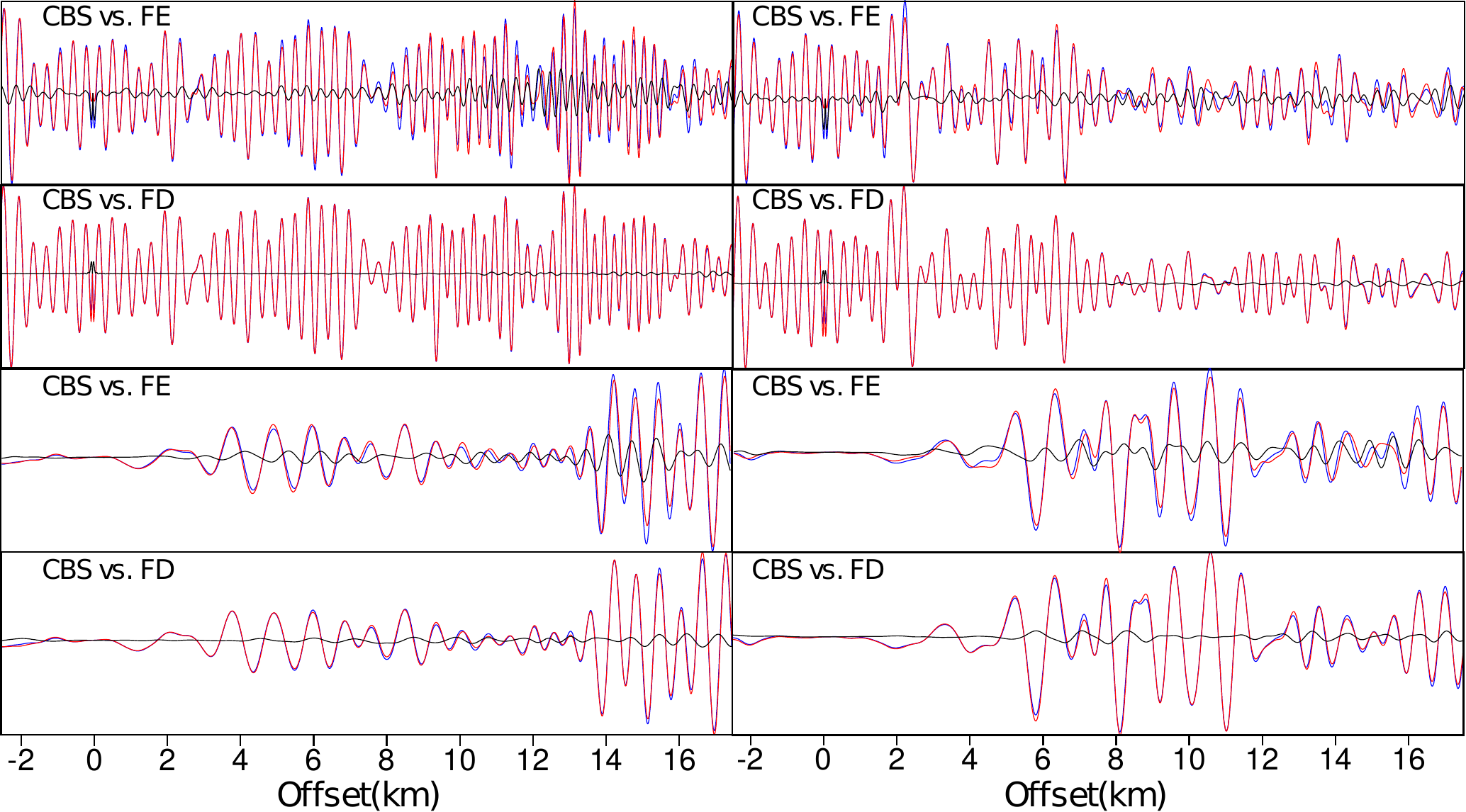}
\caption{3D EAGE/SEG overthrust model. Same as Figure~\ref{fig_over_log} but the FEFD simulation is performed on regular mesh with a discretization rule of four grid point per minimum wavelength.}
\label{fig_over_log_car}
\end{center}
\end{figure}

We now benchmark the  FEFD/FDFD ORAS solver against the heterogeneous overthrust model involving complex compressive tectonics (Table~\ref{tab_bench_spec}). 
The source is positioned at (2500~m,  2500~m, 500~m). We perform the FEFD and FDFD simulations with 2400 and 396 cores, respectively. The CBS, FEFD and FDFD wavefields are shown in Figure~\ref{fig_over_wavefield}. Direct comparison between these wavefields along selected profiles are shown in Fig. \ref{fig_over_log}. In Figure~\ref{fig_over_wavefield}, we show two depth slices of the wavefields at 500~m depth across the source and 3500~m depth just above the basement and two vertical sections at y=2500~m across the source and y=15000~m. In Figure~\ref{fig_over_log}, we show a direct comparison between the CBS and the FEFD/FDFD wavefields along four horizontal profiles at (x,z)=(2000~m, 500~m), (x,z)=(15000~m, 3500~m), (y,z)=(2500~m, 500~m), (y,z)=(15000~m, 3500~m). For this benchmark, the FDFD stencil clearly outperforms the FEFD counterpart in terms of accuracy. The scalar time of the FDFD simulation is also 14 times smaller than that of the $h$-adaptive FEFD simulation (Table~\ref{tab_bench_results}). A FEFD simulation on a regular mesh using four points per minimum wavelength confirms that the adaptive meshing doesn't impact the accuracy of the simulation (Figures~\ref{fig_over_wavefield_car} and \ref{fig_over_log_car} and Table~\ref{tab_bench_results}). The computational saving provided by the adaptive mesh is 3.8 for this benchmark (Table~\ref{tab_bench_results}).


\subsubsection{Deep crustal GO\_3D\_OBS model}
We now move to the GO\_3D\_OBS model. The source is positioned on the sea bottom at (10500~m, 12300~m, 900~m). It is reminded that the simulated frequency is 3.75~Hz (Table~\ref{tab_bench_spec}). We perform the FEFD and FDFD simulations using 2400 and 660 cores, respectively. In Figure~\ref{fig_crust_profile}, a direct comparison between the CBS and FEFD/FDFD wavefields is performed along with three horizontal profiles in the $y$ direction running across the source position and cross-cutting the accretionary wedge at 6~km depth, the subduction megathrust at 10~km depth and the Moho at 16~km depth. For this benchmark too, the FDFD solution outperforms the FEFD counterpart in terms of accuracy. However, the accuracy ratio between FEFD and FDFD (Err$_{fe}$/Err$_{fd}$) is smaller in the GO\_3D\_OBS model compared to the overthrust model (2 against 4.4). This probably result from the fact that the structures have a larger scale relatively to the wavelength in the GO\_3D\_OBS model relative to the overthrust counterpart. The computational cost of the FDFD scheme remains one order of magnitude higher than that of the FDFD counterpart.
%
%
\begin{figure}[htb!]
\begin{center}
\includegraphics[width=16cm,clip=true]{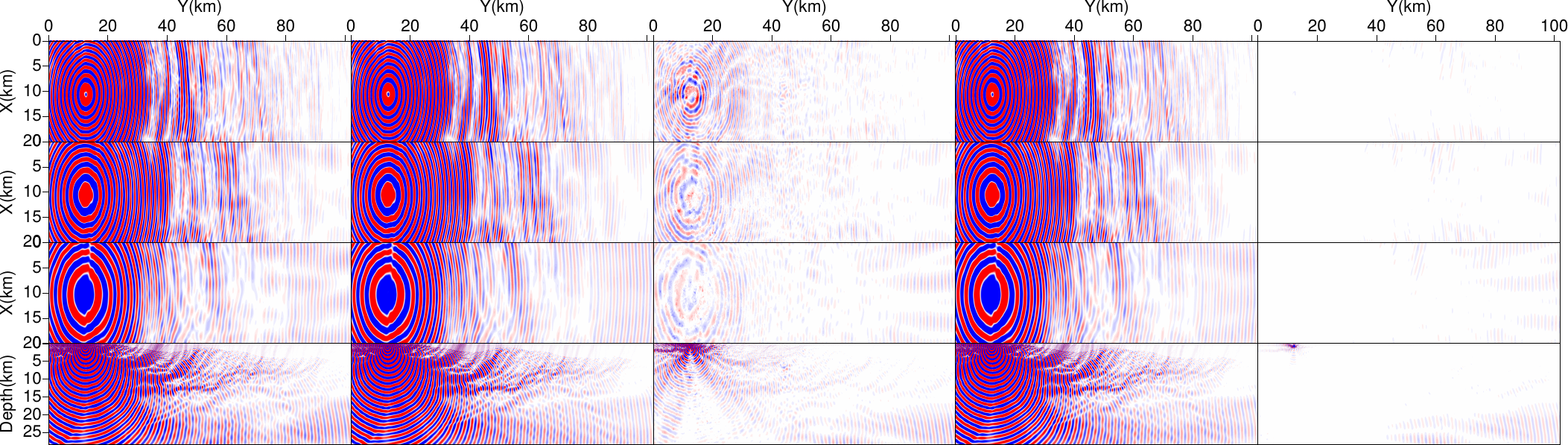}
\caption{3D GO\_3D\_OBS model. Comparison between the wavefields computed with the CBS,  FEFD and FDFD methods. The rows show from top to bottom three depth slices  at 6~km, 10~km and 16~km depth, and one vertical section at x=10~km. From left to right, the columns show the real part of the CBS wavefield, the FEFD wavefield, the differences between the two, the FDFD wavefield and its differences with the CBS wavefield.}
\label{fig_crust_wavefield}
\end{center}
\end{figure}
\begin{figure}[htb!]
\begin{center}
\includegraphics[width=13cm,clip=true]{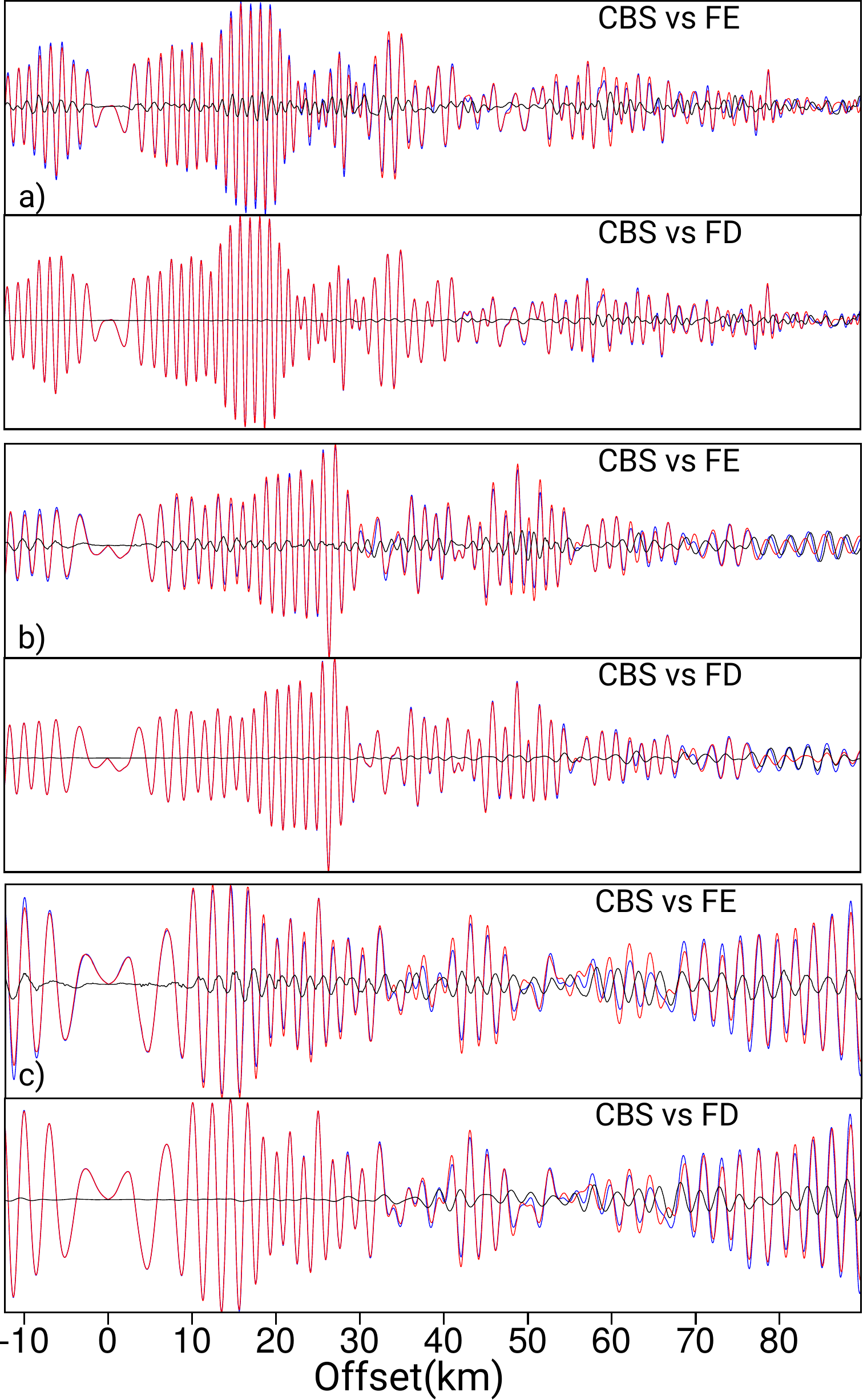}
\caption{3D GO\_3D\_OBS model. Direct comparison between the wavefields computed with the CBS (blue), FEFD/FDFD (red) methods. The differences are plotted in black. ($a_{1:2}$): Profile at (x,z)=(10k~m, 6~km). ($a_{3:4}$): Profile at (x,z)=(10~km, 10~km).  . ($a_{5:6}$): Profile at (x,z)=(10~km, 16~km). $a_{1}, a_{3}, a_{5}$: FDFD solution. $a_{2}, a_{4}, a_{6}$: FEFD solution.}
\label{fig_crust_profile}
\end{center}
\end{figure}

\subsection{Weak and strong scalability analysis}
We conclude the numerical section with a weak and strong scaling analysis of the ORAS solver. We perform these analyzes with the GO\_3D\_OBS model and the FDFD method. The finite-difference grid is discretized with four grid points per minimum wavelength. According to the conclusions of the previous section, we use single precision arithmetic, we perform local factorizations with MUMPS$_{BLR}$ using $\varepsilon_{BLR}=10^{-3}$ and CGS is used for the orthogonalization of the Krylov basis.

\subsubsection{Weak scaling}
\begin{figure}[htb!]
\begin{center}
\includegraphics[width=12cm]{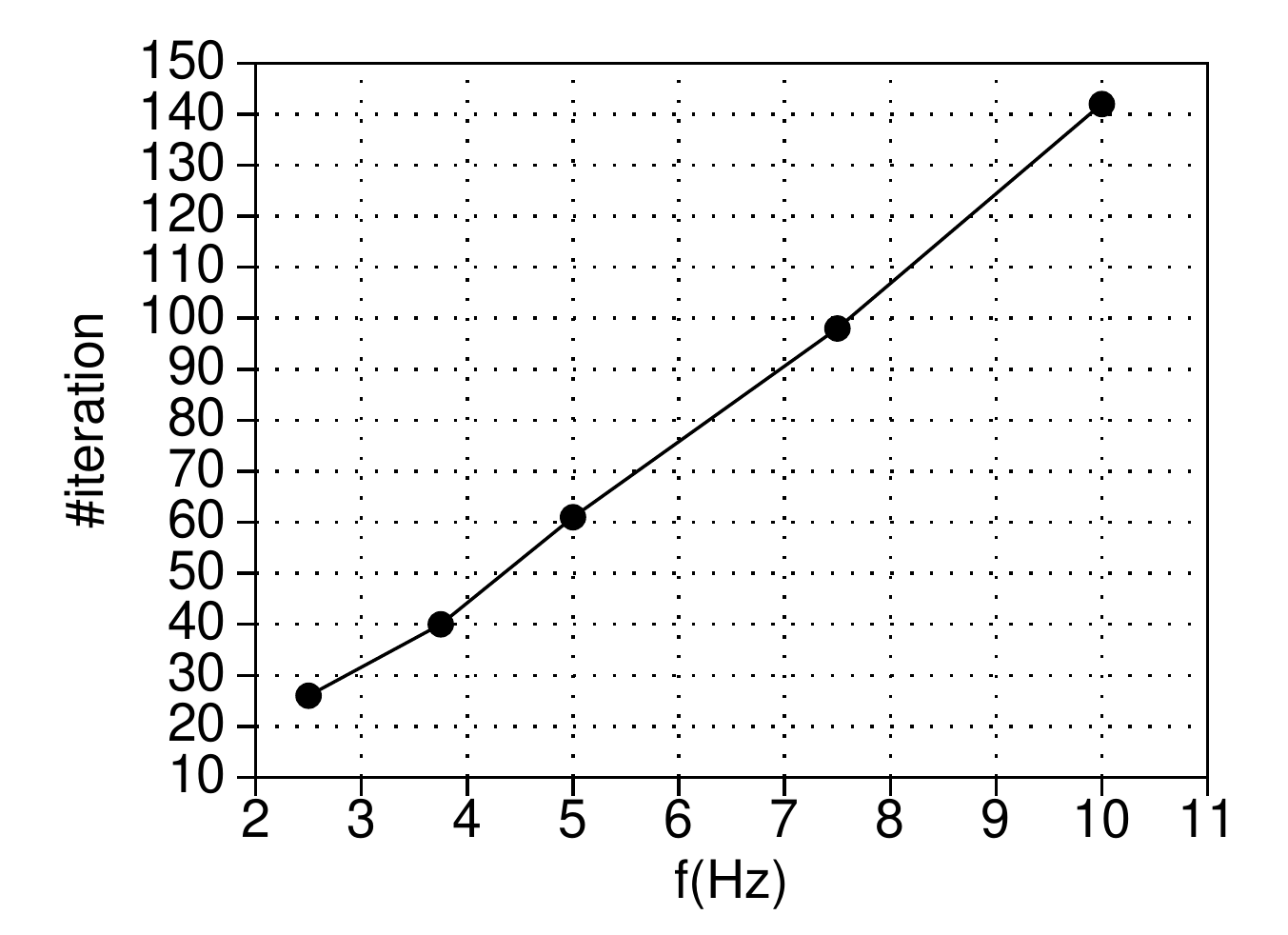}
\caption{GMRES iteration count versus frequency extracted from Table~\ref{tab_weak_scaling}. For each frequency, the number of iterations performed for the smallest number of cores is plot.}
\label{fig_itversusf}
\end{center}
\end{figure}

For the weak scaling analysis, we use frequencies ranging between 2.5~Hz and 10~Hz (Table~\ref{tab_weak_scaling}). The number of dofs (\#dof) in the Cartesian grid ranges between 21.4 and 1,160.6 million. Each simulation is performed such that the ratio \#dof/\#cores is roughly constant where \#cores denotes the number of cores. We define the weak-scaling parallel efficiency as
\begin{equation}
E_{w}=\frac{\left[ T_{tot}^{(ref)} \times \#cores^{(ref)} / \#dof^{(ref)} \right]}{ \left[T_{tot} \times \#cores / \#dof\right]},
\end{equation}
where the subscript $ref$ refers to the smallest tackled problem. \\
The results are highlighted in Table~\ref{tab_weak_scaling}. The number of iterations scales linearly with frequency, which is consistent with previous studies based upon multigrid preconditioner with shifted Laplacian \citep[e.g.,][]{Plessix_2017_CAT} (Figure~\ref{fig_itversusf}).

\begin{table}[ht!]
\begin{center}
\caption{Weak scaling. The  GO\_3D\_OBS model and the FDFD method are used to perform this analysis. Same nomenclature as that of Table~\ref{tab_solver}. $f$: Frequency in Hz. $ovl$: Size of the optimal overlap in terms of computational time.}
\label{tab_weak_scaling}
\begin{tabular}{|c|c|c|c|c|c|c|c|c|}
\hline
f(Hz) &  \#dof($10^6$)   & \#cores       &  \#it & ovl  &  $T_f(s)$  & $T_s(s)$      &  $T_{tot}(s)$   &  $E_{w}$      \\ \hline
2.5   & 21.4            &     60   & 26   &  3   & 52.5       & 24.8          &  77.3          &        1          \\ \hline
2.5   & 21.4            &    168   & 30   &  3   & 13.7       &  8.9          &  22.6          &        1.222      \\ \hline
3.75  & 67.5            &    360   & 40   &  3   & 21.9       & 18.6          &  40.5          &        1.003      \\ \hline
3.75  & 67.5            &    660   & 45   &  3   & 10.8       & 11.5          &  22.3          &        0.994      \\ \hline
3.75  & 67.5            &   1344   & 53   &  3   &  4.4       &  7.1          &  11.5          &        0.947      \\ \hline
3.75  & 67.5            &   2380   & 59   &  3   &  2.5       &  4.8          &  7.3           &        0.842      \\ \hline
5     & 153.5           &   875    & 61   &  3   & 20.9       & 26.0          &  46.9          &        0.811      \\ \hline
5     & 153.5           &  1344    & 63   &  3   & 12.4       & 18.6          &  31.0          &        0.798      \\ \hline
7.5   & 500.1           &  2450    & 98   &  4   & 26.7       & 60.8          &  87.5          &        0.506      \\ \hline
7.5   & 500.1           &  4752    & 108  &  4   & 12.2       & 39.7          &  51.9          &        0.439      \\ \hline
10    &1160.6           &  9856    & 142  &  5   & 15.1       & 70.8          &  85.9          &        0.297      \\ \hline
\end{tabular}
\end{center}
\end{table}
%
%

\subsubsection{Strong scaling}
\begin{figure}[ht!]
    \begin{center}
    \includegraphics[width=0.49\linewidth,clip=true]{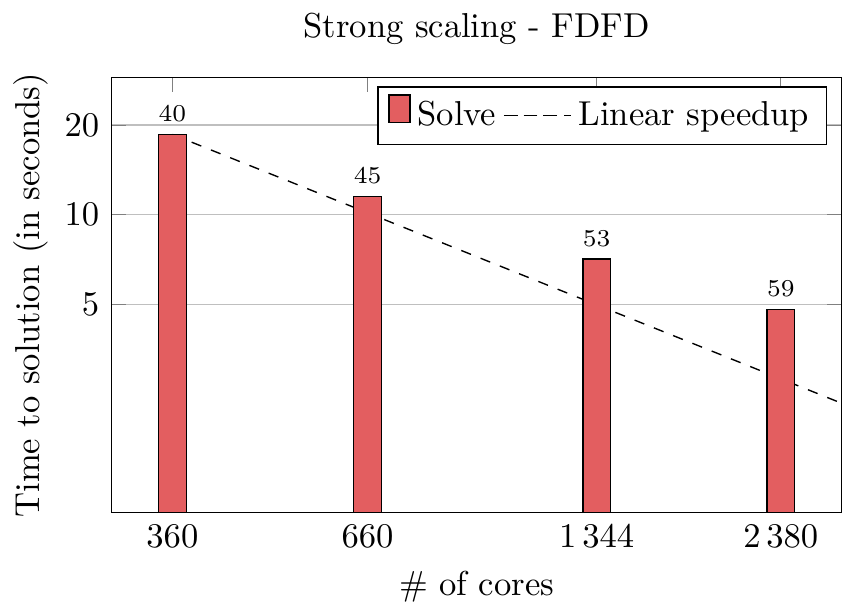}
    \includegraphics[width=0.49\linewidth,clip=true]{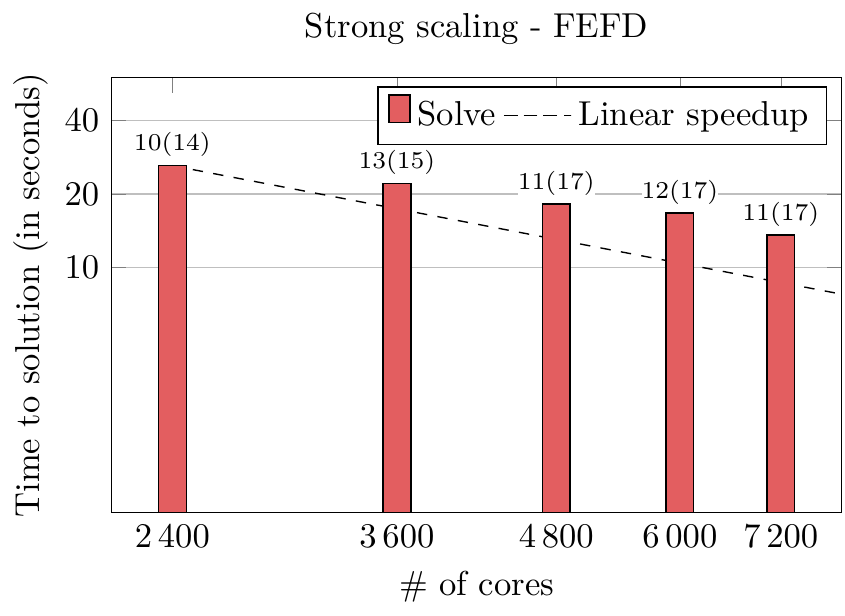}
    \caption{Strong scaling test for the FDFD (left) and FEFD (right) methods on the GO\_3D\_OBS model. Scalability of the GMRES solution time $T_s(s)$. Number of GMRES iterations is also reported, with average number of inner iterations to solve the coarse problem in parentheses for the FEFD method.}
    \label{fig_strong_scaling}
\end{center}
\end{figure}

We perform a strong scaling analysis of the two-level ORAS preconditioner in the finite element case for the 3.75~Hz frequency. The number of cores goes from 2400 to 7200 (Table~\ref{tab_strong_scaling}). \\
The parallel efficiency is defined by
\begin{equation}
E_{s}=\frac{\left[ T_{tot}^{(ref)} \times \#cores^{(ref)} \right]}{\left[T_{tot} \times \#cores \right]},
\end{equation}
The results outlined in Table~\ref{tab_strong_scaling} highlight the high scalability of the ORAS solver. Figure~\ref{fig_strong_scaling} also illustrates the scalability of the GMRES solution time for the FDFD and FEFD methods for this test case (results for the FDFD method are extracted from Table~\ref{tab_weak_scaling} for the 3.75~Hz frequency). We can see that the number of GMRES iterations for the FEFD method is more stable thanks to the coarse grid component of the preconditioner.
\begin{table}[ht!]
\begin{center}
\caption{Strong scaling. The  GO\_3D\_OBS model and the FEFD method are used to perform this analysis. Same nomenclature as that of Table~\ref{tab_weak_scaling}. $\#it$: Number of GMRES iterations, with average number of inner iterations to solve the coarse problem in parentheses.}
\label{tab_strong_scaling}
\begin{tabular}{|c|c|c|c|c|c|c|c|c|c|}
\hline
f(Hz)     &  \#dof($10^6$)    &  \#cores & \#it  & ovl  &  $T_f(s)$ & $T_s(s)$  & $T_{tot}(s)$ & $E_s$ \\ \hline
3.75      & 597.7  &  2400   & 10(14)   &  2   &   62.3    & 26.1      & 88.4   &  1      \\ \hline
3.75      & 597.7  &  3600   & 13(15)   &  2   &   38.4    & 22.1      & 60.5   &  0.974   \\ \hline
3.75      & 597.7  &  4800   & 11(17)   &  2   &   28.4    & 18.2      & 46.6   &  0.948   \\ \hline
3.75      & 597.7  &  6000   & 12(17)   &  2   &   19.4    & 16.7      & 36.1   &  0.980   \\ \hline
3.75      & 597.7  &  7200   & 11(17)   &  2   &   18.1    & 13.6      & 31.7   &  0.930   \\ \hline
\end{tabular}
\end{center}
\end{table}

\section{Discussion}
We propose a forward engine to perform 3D frequency-domain FWI of datasets collected by long-offset seabed acquisitions implemented with sparse layout of ocean bottom nodes.
In particular, we design carefully the experimental setup allowing us to achieve the best compromise between computational efficiency and simulation accuracy, the latter being quantitatively controlled by the solutions of the CBS method \citep{Osnabrugge_2016_CBS}. We also compare two discretization schemes based on the recently-proposed 27-point wavelength-adaptive finite-difference (FDFD) method on regular Cartesian grid \citep{Aghamiry_2021_AFD} and an $h$-adaptive finite-element method (FEFD) on unstructured tetrahedral mesh. 

\subsection{Which solver for FWI of sparse long-offset data?}
The pros and cons of the ORAS solver for FWI applications can be discussed relative to sparse direct solvers and explicit time marching methods. Theoretical complexities of each method for a dense 3D fixed-spread acquisition are basically identical and hence are useless to rate them in a relative sense \citep{Nihei_2007_FRM,Plessix_2007_HIS,Virieux_2009_TLE,Sourbier_2011_GEOP,Plessix_2017_CAT,Kostin_2019_DFA}.  Instead, periodic assessment of the methods carried out with realistic numerical experiments that incorporate the latest advances in the field of linear algebra may be more instructive. Such a comparative analysis is presented in \citet{Plessix_2017_CAT} where a Krylov-subspace iterative solver based upon BI-CGSTAB \citep{VanderVost_1992_BIC} and multigrid preconditioner with complex shifted-Laplacian \citep{Riyanti_2007_PMP} is assessed against an explicit time marching method.
To gain an up-to-date insight on the relative performances of the three possible approaches (namely, the two frequency-domain approaches based upon either a preconditioned iterative solver or a pure direct solver and the time marching method), we compare the results obtained with the ORAS solver for the (acoustic) GO\_3D\_OBS model with 130 RHSs with those obtained with MUMPS$_{BLR}$ and an explicit time-stepping finite-difference time-domain solver based on a classical $\mathcal{O}(\Delta t^2,\Delta t^8)$ staggered-grid stencil using a grid interval of 100~m. In all cases, we use the Occigen supercomputer. \\

For the FDFD simulation with MUMPS$_{BLR}$, we use $\varepsilon_{BLR}=10^{-5}$, which is required to reach the desired accuracy \citep{Amestoy_2021_EDA}. Note that a relaxed value of $\varepsilon_{BLR}=10^{-3}$ was used to solve the local problems with ORAS since a high accuracy is not required to build the preconditioner. We perform the simulation with MUMPS$_{BLR}$ on 80 nodes using one MPI process per node and 24 threads per process to minimize memory overhead while exploiting multithreading for computational efficiency. The elapsed time to perform LU factorization and compute 130 solutions is 1000~s and 28.3~s, respectively, leading to a total elapsed time of 1028.3s and a scalar time of 548 scalar hours. \\
The FDTD simulation relies on  two-level MPI parallelism with an embarrassing parallelism implemented by distributing the sources over the processes and a domain decomposition of the computational domain. We use a total of 12,480 computer cores to implement this two-level parallelism. Two sequential paths over the RHSs are performed to process 65 RHSs in parallel while the domain is decomposed into  $2 \times 24 \times 4 = 192$ subdomains. The parallel runtime to perform one simulation is 132s. Hence, the parallel run time to perform 130 simulations is 264~s and the scalar time is 915~hours. These numbers can be compared with those obtained with ORAS$_{FDFD}$ with a scalar time of 98 core hours. This number is inferred from the last line of Table~\ref{tab_solver} by multiplying the parallel run time for 130 RHSs ($T_{tot} = 261.7~s$) by the number of cores (1344).

Before comparing these numbers, it is worth reminding that, during a FWI iteration, a time-domain method generally requires to recompute the incident wavefields backward in time from the boundaries during the adjoint simulation unlike frequency-domain methods where the incident and adjoint wavefields associated with a block of RHSs can be kept in memory. Alternatively, check-pointing strategies can be used to perform this re-computation, in particular in attenuating media where time-reverse modeling of the incident wavefields are unstable \citep{Symes_2007_RTM,Anderson_2012_TCM,Yang_2016_CAR}. Moreover, attenuation significantly increases the cost of FDTD simulation (by a factor of around 3 according to Figure 13 of \citet{Plessix_2017_CAT}), while physical attenuation will improve the convergence rate of the ORAS solver. As an illustrative example, when a 3.75-Hz simulation is performed in the GO\_3D\_OBS model using the $Q_P$ model shown in the Figure~6d of \citet{Gorszczyk_2021_GNT}, the iteration count drops from 45 to 39 and the solution step for 1 RHS from 11.5~s to 9.8~s, hence leading to a computational saving of around 15\%. However, it is worth reminding that small patches of frequencies (typically, two to four) are generally inverted during frequency-domain FWI in particular when attenuation reconstruction is targeted \citep[e.g.][]{Operto_2018_MFF}. Therefore, the scalar time of frequency-domain modelling should be multiplied by the number of frequencies involved in the patch (these frequencies are ideally processed in parallel with an embarassing parallelism) for a fair comparison with time-domain methods, the forward engine of which provides the full bandwidth of the wavefield in one go. 
Concerning the two frequency-domain approaches, we note that the speedup achieved during the block processing of multiple RHSs is 18 and 3.6 for MUMPS and ORAS, respectively. Indeed, direct solvers are more suitable for dense acquisitions to really exploit their ability to process a large number of sparse RHSs. 

Lumping all these specifications together, we conclude that the ORAS-FDFD solver should be highly competitive relative to a time-domain solver in particular when attenuation effects should be taken into account. It allows for compact volume of data to be managed, attenuation speeds up the convergence rate of the solver and the good strong scalability provides the necessary flexibility to optimally tune the two-level parallelism based on frequency distribution and domain decomposition according to the acquisition specifications, target size and available computational resources. The solver MUMS$_{BLR}$ further contributes to mitigate the memory cost of the method and the computational efficiency of the local subproblems. Overall, the amount of resources used to perform the simulations (\#cores in Table~\ref{tab_comparison}) shows that frequency-domain approaches are more suitable when a more limited amount of resources are available, which is consistent with the conclusions of \citet{Kostin_2019_DFA}.
\begin{table}[ht!]
\begin{center}
\caption{Cost of the FDTD, MUMPS-FDFD and ORAS-FDFD simulations in the GO\_3D\_OBS model for a frequency of 3.75~Hz. $\#cores$: Number of computer cores. $T_s^{(130)}(s)$: Elapsed time to perform 130 simulations in parallel. $T_{hc}(hr)$: scalar hours (elapsed time times the number of cores). $T_f(s)$: Elapsed time to perform parallel LU factorization with MUMPS-FDFD. $T_s^1(s)$: Elapsed time to process one RHS with MUMPS-FDFD.}
\label{tab_comparison}
\begin{tabular}{|c|c|c|c|c|c|c|c|c|c|c|}
\hline
\multicolumn{3}{|c|}{FDTD} & \multicolumn{5}{|c|}{MUMPS-FDFD} & \multicolumn{3}{|c|}{ORAS-FDFD} \\ \hline
\#c   &  $T_s^{(130)}(s)$ & $T_{hc}(hr)$ & \#cores & $T_f(s)$ & $T_s^1(s)$ & $T_s^{130}(s)$ & $T_{hc}(hr)$ & \#cores & $T_s^{130}(s)$ & $T_{hc}(hr)$ \\ \hline
12480 &  264  & 915 & 1920  & 1000 & 3.9 &  28.3 & 548 & 1344 &  256.8 & 96 \\ \hline
\end{tabular}
\end{center}
\end{table}

\subsection{Finite elements versus finite differences}
Our results suggest that the FDFD method outperforms the FEFD counterpart in terms of computational efficiency and accuracy for the four acoustic benchmarks tackled in the previous section. This results because the recently proposed wavelength-adaptive 27-point FDFD stencil \citep{Aghamiry_2021_AFD} leads to a quite significant improvement of accuracy compared to the non-adaptive version \citep{Operto_2007_FDFD} without introducing computational overhead.
In the frame of FWI, a clear advantage of the FEFD method is however the flexibility to comply the mesh with known boundary such as bathymetry. Generally, FWI is performed with a frequency continuation strategy to mitigate cycle skipping. The grid interval can be matched to the frequency to regularize the inversion and mitigates the computational cost \citep{Operto_2015_ETF}. However, when FDFD method on regular grid is used, this grid adaptation may lead to inaccurate representation of the bathymetry when low frequencies are processed on coarse grids. This can inject shallow inaccuracies in the velocity model near the seabed, which are propagated deeper in the velocity model and which are difficult to correct at higher frequencies, hence making the inversion stall at high frequencies. The FEFD method will clearly give the necessary flexibility to comply the mesh with the bathymetry whatever the size of the elements, hence facilitating the implementation of multiscale approaches. Moreover, elastic physics will further increase the value of the $h$-adaptive meshing used by FEFD methods to account for the low shear wavespeeds in the near surface. A re-assessment of the two discretization schemes is clearly necessary for elastic physics.


\section{Conclusion}
We have proposed an efficient, accurate and versatile forward engine for 3D frequency-domain FWI from ultra-long offset stationary-recording survey. The solver relies on an iterative solver and a domain-decomposition based preconditioner where the local subproblems are tackled with sparse direct solver. Two different discretization schemes based on wavelength-adaptive finite differences on regular Cartesian grid and $h$-adaptive finite elements on unstructured tetrahedral meshes can be interfaced with the solver. We perform a careful parametric analysis that allows us to define the experimental setup providing the best compromise between computational efficiency and accuracy. The strong scalability of the solver provides the necessary flexibility to implement two levels of parallelism by frequency distribution and domain decomposition according to the specifications of the acquisition. The perspective of this work is to use this solver as a forward engine for FWI to perform 3D visco-acoustic case studies from sparse long-offset node acquisitions. Extension to visco-elastic physics is another objective. 


\section*{Acknowledgments}
This study was granted access to the HPC resources of SIGAMM (http://crimson.oca.eu) and CINES/IDRIS under the allocation 0596 made by GENCI. This study was partially funded by the WIND consortium (https://www.geoazur.fr/WIND) sponsored by Chevron, Shell and Total. The finite-difference time-domain simulation of the discussion section has been performed by L. Combe (Geoazur) with the GeoInv3D\_fwk software. We thank the MUMPS team (P. Amestoy (Mumps Technologies, ENS Lyon), A. Buttari (CNRS, IRIT), J.-Y. L'Excellent (Mumps Technologies, ENS Lyon), T. Mary (Sorbonne Universit\'e, CNRS, LIP6) and C. Puglisi (Mumps Technologies, ENS Lyon)) for proving the MUMPS solver (https://mumps-consortium.org) and for fruitful discussions.


\bibliographystyle{seg}  

\begin{thebibliography}{}
\itemsep0pt

\bibitem[Aghamiry et~al., 2021]{Aghamiry_2021_AFD}
Aghamiry, H., A. Gholami, L. Combe, and S. Operto,  2021, {Accurate 3D
  frequency-domain seismic wave modeling with the wavelength-adaptive 27-point
  finite-difference stencil: a tool for full waveform inversion}: Geophysics,
  {\bf Submitted, 	arXiv:2108.08730}.

\bibitem[Ainsworth and Wajid, 2010]{Ainsworth:2010:OBS}
Ainsworth, M., and H. Wajid,  2010, Optimally blended spectral-finite element
  scheme for wave propagation and nonstandard reduced integration: SIAM Journal
  on Numerical Analysis, {\bf 48}, 346--371.

\bibitem[Alappat et~al., 2020]{Alappat_2020_RAC}
Alappat, C., G. Hager, O. Schenk, J. Thies, A. Basermann, A. Bischop, H.
  Fehske, and G. Wellein,  2020, A recursive algebraic coloring technique for
  hardware-efficient symmetric sparse matrix-vector multiplication: ACM
  Transactions on Parallel Computing, {\bf 7}, 1--37.

\bibitem[Amestoy et~al., 2016]{Amestoy_2016_FFF}
Amestoy, P., R. Brossier, A. Buttari, J.-Y. L'Excellent, T. Mary, L.
  M\'etivier, A. Miniussi, and S. Operto,  2016, Fast {3D} frequency-domain
  {FWI} with a parallel {B}lock {L}ow-{R}ank multifrontal direct solver:
  application to {OBC} data from the {N}orth {S}ea: Geophysics, {\bf 81}, R363
  -- R383.

\bibitem[Amestoy et~al., 2018]{Amestoy_2018_ESM}
Amestoy, P., A. Buttari, J.~Y. L'Excellent, and T. Mary,  2018, On exploiting
  sparsity of multiple right-hand sides in sparse direct solvers: {SIAM}
  {J}ournal on {S}cientific {C}omputing, {\bf 41}, A269--A291.

\bibitem[Amestoy et~al., 2021]{Amestoy_2021_EDA}
Amestoy, P.~R., A. Buttari, L. Combe, M. Gerest, J.-Y. {L'Excellent}, T. Mary,
  S. Operto, and C. Puglisi,  2021, {Up to date assessment of 3D
  frequency-domain full waveform inversion based on the sparse multifrontal
  solver MUMPS}: Presented at the {Fifth EAGE Workshop on High Performance
  Computing for Upstream}.

\bibitem[Amestoy et~al., 2019a]{Amestoy_2019_PSB}
Amestoy, P.~R., A. Buttari, J.-Y. {L'Excellent}, and T. Mary,  2019a,
  {Performance and Scalability of the Block Low-Rank Multifrontal Factorization
  on Multicore Architectures}: ACM Transactions on Mathematical Software, {\bf
  45}, 2:1--2:26.

\bibitem[Amestoy et~al., 2019b]{Amestoy_2019_EUS}
Amestoy, P.~R., S. {de la Kethulle de Ryhov}, J.-Y. {L'Excellent}, G. Moreau,
  and D.~V. Shantsev,  2019b, {Efficient use of sparsity by direct solvers
  applied to 3D controlled-source EM problems}: Computational Geosciences, {\bf
  23}, 1237--1258.

\bibitem[Aminzadeh et~al., 1997]{Aminzadeh_1997_DSO}
Aminzadeh, F., J. Brac, and T. Kunz,  1997, {3-D} {Salt} and {Overthrust}
  models: {SEG/EAGE} 3-{D} {M}odeling {S}eries {N}o.1.

\bibitem[Anderson et~al., 2012]{Anderson_2012_TCM}
Anderson, J.~E., L. Tan, and D. Wang,  2012, Time-reversal checkpointing
  methods for {RTM} and {FWI}: Geophysics, {\bf 77}, S93--S103.

\bibitem[Babuska and Sauter, 1997]{Babuska:1997:IPE}
Babuska, I.~M., and S.~A. Sauter,  1997, Is the pollution effect of the {FEM}
  avoidable for the {H}elmholtz equation considering high wave numbers?: SIAM
  Journal on Numerical Analysis, {\bf 34}, 2392--2423.

\bibitem[{Ben Hadj Ali} et~al., 2011]{Ben-Hadj-Ali_2011_GEO}
{Ben Hadj Ali}, H., S. Operto, and J. Virieux,  2011, An efficient
  frequency-domain full waveform inversion method using simultaneous encoded
  sources: Geophysics, {\bf 76}, R109.

\bibitem[B\'erenger, 1994]{Berenger_1994_PML}
B\'erenger, J.-P.,  1994, A perfectly matched layer for absorption of
  electromagnetic waves: Journal of Computational Physics, {\bf 114}, 185--200.

\bibitem[B\"ollhofer et~al., 2019]{Bollofer_2019_LSI}
B\"ollhofer, M., A. Eftekhari, S. Scheidegger, and O. Schenk,  2019,
  Large-scale sparse inverse covariance matrix estimation: SIAM J. Sci.
  Comput., {\bf 41(1)}, A380–A401.

\bibitem[B\"ollhofer et~al., 2020]{Bollofer_2020_SSD}
B\"ollhofer, M., O. Schenk, R. Janalik, S. Hamm, and K. Gullapalli,  2020,
  State-of-the-art sparse direct solvers: Parallel Algorithms in Computational
  Science and Engineering, Birkhäuser,  3--33.

\bibitem[Bonazzoli et~al., 2019]{Bonazzoli:2019:ADD}
Bonazzoli, M., V. Dolean, I.~G. Graham, E.~A. Spence, and P.-H. Tournier,
  2019, A 2-level domain decomposition preconditioner for the time-harmonic
  {Maxwell}'s equations: Math. Comp., {\bf 88}, 2559--2604.

\bibitem[Bootland and Dolean, 2019]{Bootland:2019:ODN}
Bootland, N., and V. Dolean,  2019, On the {D}irichlet-to-{N}eumann coarse
  space for solving the {H}elmholtz problem using domain decomposition:
  Presented at the Numerical Mathematics and Advanced Applications ENUMATH
  2019.

\bibitem[Bootland et~al., 2021a]{Bootland:2021:SWA}
Bootland, N., V. Dolean, P. Jolivet, F. Nataf, S. Operto, and P.-H. Tournier,
  2021a, Several ways to achieve robustness when solving wave propagation
  problems.

\bibitem[Bootland et~al., 2021b]{Bootland:2021:ACS}
Bootland, N., V. Dolean, P. Jolivet, and P.-H. Tournier,  2021b, A comparison
  of coarse spaces for {H}elmholtz problems in the high frequency regime:
  Computers and Mathematics with Applications, {\bf 98}, 239--253.

\bibitem[Brossier et~al., 2010]{Brossier_2010_FNM}
Brossier, R., V. Etienne, S. Operto, and J. Virieux,  2010, Frequency-domain
  numerical modelling of visco-acoustic waves based on finite-difference and
  finite-element discontinuous galerkin methods, {\it in} Acoustic Waves:
  SCIYO,  125--158.

\bibitem[Carcione and Robinson, 2002]{Carcione_2002_AEA}
Carcione, J.~M., and E.~M. Robinson,  2002, On the acoustic-electromagnetic
  analogy for the reflection-refraction problem: Studia Geophysica et
  Geodaetica, {\bf 46}, 321--346.

\bibitem[Chen et~al., 2012]{Chen_2012_DMF}
Chen, Z., D. Cheng, and T. Wu,  2012, {A dispersion minimizing finite
  difference scheme and preconditioned solver for the 3D Helmholtz equation}:
  Journal of Computational Physics, {\bf 231}, 8152--8175.

\bibitem[Chew and Weedon, 1994]{Chew_1994_PMM}
Chew, W.~C., and W.~H. Weedon,  1994, A 3-{D} perfectly matched medium from
  modified {M}axwell's equations with stretched coordinates: Microwave and
  Optical Technology Letters, {\bf 7}, 599--604.

\bibitem[Conen et~al., 2014]{Conen:2014:ACS}
Conen, L., V. Dolean, R. Krause, and F. Nataf,  2014, A coarse space for
  heterogeneous {Helmholtz} problems based on the {Dirichlet-to-Neumann}
  operator: J. Comput. Appl. Math., {\bf 271}, 83--99.

\bibitem[Dapogny et~al., 2014]{dapogny2014three}
Dapogny, C., C. Dobrzynski, and P. Frey,  2014, Three-dimensional adaptive
  domain remeshing, implicit domain meshing, and applications to free and
  moving boundary problems: Journal of computational physics, {\bf 262},
  358--378.

\bibitem[Dolean et~al., 2015]{Dolean:2015:IDD}
Dolean, V., P. Jolivet, and F. Nataf,  2015, An introduction to domain
  decomposition methods: Society for Industrial and Applied Mathematics (SIAM),
  Philadelphia, PA.
\newblock (Algorithms, theory, and parallel implementation).

\bibitem[Dolean et~al., 2020a]{Dolean:2020:IFD}
Dolean, V., P. Jolivet, P.-H. Tournier, and S. Operto,  2020a, Iterative
  frequency-domain seismic wave solvers based on multi-level
  domain-decomposition preconditioners: Presented at the 82$^{th}$ Annual EAGE
  Meeting (Amsterdam).

\bibitem[Dolean et~al., 2020b]{Dolean:2020:LSF}
--------, 2020b, Large-scale frequency-domain seismic wave modeling on {\it
  {h}}-adaptive tetrahedral meshes with iterative solver and multi-level
  domain-decomposition preconditioners: Presented at the SEG2020 Annual
  Meeting.
\newblock (Houston).

\bibitem[Duff et~al., 2017]{Duff:2017:DMS}
Duff, I.~S., A.~M. Erisman, and J.~K. Reid,  2017, Direct methods for sparse
  matrices, second edition: Oxford University Press, Inc.

\bibitem[Erlangga, 2008]{Erlangga:2008:AIM}
Erlangga, Y.,  2008, Advances in iterative methods and preconditioners for the
  {Helmholtz} equation: Archives of Computational Methods in Engineering, {\bf
  15}, 37--66.

\bibitem[Ernst and Gander, 2012]{Ernst:12:NAM}
Ernst, O.~G., and M.~J. Gander,  2012, Why it is difficult to solve {H}elmholtz
  problems with classical iterative methods, {\it in} Numerical analysis of
  multiscale problems: Springer,  325--363.

\bibitem[Gander and Nataf, 2000]{Gander:2000:AILU}
Gander, M., and F. Nataf,  2000, {AILU} for {H}elmholtz problems: a new
  preconditioner based on an analytic factorization: Comptes Rendus de
  l'Académie des Sciences - Series I - Mathematics, {\bf 331}, 261--266.

\bibitem[Gander and Zhang, 2019]{Gander:2018:SIREV}
Gander, M.~J., and H. Zhang,  2019, A class of iterative solvers for the
  {H}elmholtz equation: Factorizations, sweeping preconditioners, source
  transfer, single layer potentials, polarized traces, and optimized {S}chwarz
  methods: SIAM Review, {\bf 61}, 3--76.

\bibitem[Giraud et~al., 2021]{Giraud_2021_BGC}
Giraud, L., Y.-F. Jing, and Y. Xiang,  2021, A block minimum residual norm
  subspace solver for sequences of multiple left and right-hand side linear
  systems: Technical report, RR-9393, Inria Bordeaux Sud-Ouest, hal-03146213v2.

\bibitem[G\'orszczyk and Operto, 2021]{Gorszczyk_2021_GNT}
G\'orszczyk, A., and S. Operto,  2021, {GO\_3D\_OBS}: the multi-parameter
  benchmark geomodel for seismic imaging method assessment and next-generation
  3d survey design (version 1.0): Geoscientific Model Development, {\bf 14},
  1773–1799.

\bibitem[G{\'{o}}rszczyk et~al., 2017]{Gorszczyk_2017_TRW}
G{\'{o}}rszczyk, A., S. Operto, and M. Malinowski,  2017, Toward a robust
  workflow for deep crustal imaging by {FWI} of {OBS} data: The eastern nankai
  trough revisited: Journal of Geophysical Research: Solid Earth, {\bf 122},
  4601--4630.

\bibitem[Gosselin-Cliche and Giroux, 2014]{Gosselin_2014_FDF}
Gosselin-Cliche, B., and B. Giroux,  2014, {3D frequency-domain
  finite-difference viscoelastic-wave modeling using weighted average 27-point
  operators with optimal coefficients}: Geophysics, {\bf 79}, T169--T188.

\bibitem[Graham et~al., 2017]{Graham:2017:RRD}
Graham, I.~G., E.~A. Spence, and E. Vainikko,  2017, Recent results on domain
  decomposition preconditioning for the high-frequency helmholtz equation using
  absorption: Lahaye D., Tang J., Vuik K. (eds) Modern Solvers for Helmholtz
  Problems. Geosystems Mathematics. Birkhäuser, Cham,  3--26.

\bibitem[Jolivet and Tournier, 2016]{Jolivet_2016_BIM}
Jolivet, P., and P.~H. Tournier,  2016, Block iterative methods and recycling
  for improved scalability of linear solvers: {SC16}: International Conference
  for High Performance Computing, Networking, Storage and Analysis, 190--203.

\bibitem[Karypis, 2013]{Karypis_2013_MSP}
Karypis, G.,  2013, {METIS - A} software package for partitioning unstructured
  graphs, partitioning meshes, and computing fill-reducing orderings of sparse
  matrices - version 5.1.0.
\newblock University of Minnesota.

\bibitem[Kostin et~al., 2019]{Kostin_2019_DFA}
Kostin, V., S. Solovyev, A. Bakulin, and M. Dmitriev,  2019, {Direct
  frequency-domain {3D} acoustic solver with intermediate data compression
  benchmarked against time-domain modeling for full-waveform inversion
  applications}: Geophysics, {\bf 84(4)}, T207--T219.

\bibitem[Kuvshinov and Mulder, 2006]{Kuvshinov_2006_EST}
Kuvshinov, B.~N., and W.~A. Mulder,  2006, The exact solution of the
  time-harmonic wave equation for a linear profile: Geophysical Journal
  International, {\bf 167}, 659--662.

\bibitem[Li et~al., 2020]{Li_2020_FEW}
Li, Y., R. Brossier, and L. M{\'{e}}tivier,  2020, {3D} frequency-domain
  elastic wave modeling with the spectral element method using a massively
  parallel direct solver: GEOPHYSICS, {\bf 85}, T71--T88.

\bibitem[Li et~al., 2015]{Li_2015_FDE}
Li, Y., L. M\'etivier, R. Brossier, B. Han, and J. Virieux,  2015, {2D and 3D}
  frequency-domain elastic wave modeling in complex media with a parallel
  iterative solver: Geophysics, {\bf 80(3)}, T101--T118.

\bibitem[Marfurt, 1984]{Marfurt_1984_AFF}
Marfurt, K.,  1984, Accuracy of finite-difference and finite-element modeling
  of the scalar and elastic wave equations: Geophysics, {\bf 49}, 533--549.

\bibitem[{MUMPS team}, 2021]{MUMPS_2021_MMP}
{MUMPS team},  2021, {MU}ltifrontal {M}assively {P}arallel {S}olver ({MUMPS}
  5.4.1) {U}sers' guide ({August}, 2021).
\newblock {Mumps Technologies}, http://mumps-solver.org.

\bibitem[Nihei and Li, 2007]{Nihei_2007_FRM}
Nihei, K.~T., and X. Li,  2007, Frequency response modelling of seismic waves
  using finite difference time domain with phase sensitive detection
  ({TD-PSD}): Geophysical Journal International, {\bf 169}, 1069--1078.

\bibitem[Operto et~al., 2014]{Operto_2014_FAT}
Operto, S., R. Brossier, L. Combe, L. M\'etivier, A. Ribodetti, and J. Virieux,
   2014, Computationally-efficient three-dimensional visco-acoustic
  finite-difference frequency-domain seismic modeling in vertical transversely
  isotropic media with sparse direct solver: Geophysics, {\bf 79(5)},
  T257--T275.

\bibitem[Operto and Miniussi, 2018]{Operto_2018_MFF}
Operto, S., and A. Miniussi,  2018, On the role of density and attenuation in
  {3D} multi-parameter visco-acoustic {VTI} frequency-domain {FWI}: an {OBC}
  case study from the {North Sea}: Geophysical Journal International, {\bf
  213}, 2037--2059.

\bibitem[Operto et~al., 2015]{Operto_2015_ETF}
Operto, S., A. Miniussi, R. Brossier, L. Combe, L. M\'etivier, V. Monteiller,
  A. Ribodetti, and J. Virieux,  2015, Efficient {3-D} frequency-domain
  mono-parameter full-waveform inversion of ocean-bottom cable data:
  application to {V}alhall in the visco-acoustic vertical transverse isotropic
  approximation: Geophysical Journal International, {\bf 202}, 1362--1391.

\bibitem[Operto et~al., 2007]{Operto_2007_FDFD}
Operto, S., J. Virieux, P. Amestoy, J.-Y. L'\'Excellent, L. Giraud, and H. {Ben
  Hadj Ali},  2007, 3{D} finite-difference frequency-domain modeling of
  visco-acoustic wave propagation using a massively parallel direct solver: A
  feasibility study: Geophysics, {\bf 72}, SM195--SM211.

\bibitem[Osnabrugge et~al., 2016]{Osnabrugge_2016_CBS}
Osnabrugge, G., S. Leedumrongwatthanakun, and I.~M. Vellekoop,  2016, A
  convergent born series for solving the inhomogeneous helmholtz equation in
  arbitrarily large media: Journal of computational physics, {\bf 322},
  113--124.

\bibitem[Parks et~al., 2006]{Parks_2006_RKS}
Parks, M., E. de~Sturler, G. Mackey, D. Johnson, and S. Maiti,  2006, Recycling
  krylov subspaces for sequences of linear systems: {SIAM} Journal of
  Scientific Computing, {\bf 28}, 1651--1674.

\bibitem[Pellegrini, 2018]{Francois_2018_SCOTCH}
Pellegrini, F.,  2018, {PT-SCOTCH} and lib{PTSCOTCH} 6.0 user's guide - version
  6.0.5.
\newblock Universit\'{e} Bordeaux.

\bibitem[Plessix, 2007]{Plessix_2007_HIS}
Plessix, R.~E.,  2007, A {H}elmholtz iterative solver for {3D} seismic-imaging
  problems: Geophysics, {\bf 72}, SM185--SM194.

\bibitem[Plessix, 2017]{Plessix_2017_CAT}
--------, 2017, Some computational aspects of the time and frequency domain
  formulations of seismic waveform inversion, {\it in} Modern solvers for
  {H}elmholtz problems, Geosystems Mathematics: Springer,  159--187.

\bibitem[Pratt et~al., 1996]{Pratt_1996_TDV}
Pratt, R.~G., Z.~M. Song, P.~R. Williamson, and M. Warner,  1996,
  Two-dimensional velocity models from wide-angle seismic data by wavefield
  inversion: Geophysical Journal International, {\bf 124}, 323--340.

\bibitem[Riyanti et~al., 2007]{Riyanti_2007_PMP}
Riyanti, C.~D., A. Kononov, Y.~A. Erlangga, C. Vuik, C. Oosterlee, R.~E.
  Plessix, and W.~A. Mulder,  2007, A parallel multigrid-based preconditioner
  for the {3D} heterogeneous high-frequency {H}elmholtz equation: Journal of
  Computational physics, {\bf 224}, 431--448.

\bibitem[Saad, 1993]{saad1993flexible}
Saad, Y.,  1993, A flexible inner-outer preconditioned gmres algorithm: SIAM
  Journal on Scientific Computing, {\bf 14}, 461--469.

\bibitem[Saad, 2003a]{Saad:2003:IMS}
--------, 2003a, Iterative methods for sparse linear systems, 2nd ed.: Society
  for Industrial and Applied Mathematics.

\bibitem[Saad, 2003b]{Saad_2003_IMS}
--------, 2003b, Iterative {M}ethods for {S}parse {L}inear {S}ystems: SIAM.

\bibitem[Sirgue and Pratt, 2004]{Sirgue_2004_EWI}
Sirgue, L., and R.~G. Pratt,  2004, Efficient waveform inversion and imaging :
  a strategy for selecting temporal frequencies: Geophysics, {\bf 69},
  231--248.

\bibitem[Sourbier et~al., 2011]{Sourbier_2011_GEOP}
Sourbier, F., A. Haiddar, L. Giraud, H. Ben-Hadj-Ali, S. Operto, and J.
  Virieux,  2011, Three-dimensional parallel frequency-domain visco-acoustic
  wave modelling based on a hybrid direct/iterative solver: Geophysical
  Prospecting, {\bf 59}, 834--856.

\bibitem[Symes, 2007]{Symes_2007_RTM}
Symes, W.~W.,  2007, Reverse time migration with optimal checkpointing:
  Geophysics, {\bf 72}, SM213--SM221.

\bibitem[Tournier et~al., 2019]{Tournier:2019:MTI}
Tournier, P.-H., I. Aliferis, M. Bonazzoli, M. de~Buhan, M. Darbas, V. Dolean,
  F. Hecht, P. Jolivet, I.~E. Kanfoud, C. Migliaccio, F. Nataf, C. Pichot, and
  S. Semenov,  2019, Microwave tomographic imaging of cerebrovascular accidents
  by using high-performance computing: Parallel Computing, {\bf 85}, 88 -- 97.

\bibitem[{Van der Vorst}, 1992]{VanderVost_1992_BIC}
{Van der Vorst}, H.~A.,  1992, {BI-CGSTAB}: a fast and smoothly converging
  variant of bi-{CG} for the solution of nonsymmetric linear systems: {SIAM}
  Journal on Scientific and Statistical Computing, {\bf 13}, 631--644.

\bibitem[Virieux and Operto, 2009]{Virieux_2009_OFW}
Virieux, J., and S. Operto,  2009, An overview of full waveform inversion in
  exploration geophysics: Geophysics, {\bf 74}, WCC1--WCC26.

\bibitem[Virieux et~al., 2009]{Virieux_2009_TLE}
Virieux, J., S. Operto, H. {Ben Hadj Ali}, R. Brossier, V. Etienne, F.
  Sourbier, L. Giraud, and A. Haidar,  2009, Seismic wave modeling for seismic
  imaging: The Leading Edge, {\bf 28}, 538--544.

\bibitem[Yang et~al., 2016]{Yang_2016_CAR}
Yang, P., R. Brossier, L. M\'etivier, and J. Virieux,  2016, Wavefield
  reconstruction in attenuating media: A checkpointing-assisted reverse-forward
  simulation method: Geophysics, {\bf 81}, R349--R362.

\bibitem[Zienkiewicz et~al., 2005]{Zienkiewicz_2005_FEM}
Zienkiewicz, O.~C., R.~L. {Taylor}, and J.~Z. Zhu,  2005, The {F}inite
  {E}lement {M}ethod: Its basis and fundamentals: Elsevier.
\newblock (6\textsuperscript{th} edition).

\end{thebibliography}

\newcommand{\SortNoop}[1]{}

\end{document}